\documentclass[11pt]{article}
\usepackage{amssymb,amsmath,latexsym}

\begin{document}

\begin{doublespace}

\newtheorem{thm}{Theorem}[section]
\newtheorem{lemma}[thm]{Lemma}

\newtheorem{defn}[thm]{Definition}
\newtheorem{prop}[thm]{Proposition}
\newtheorem{corollary}[thm]{Corollary}
\newtheorem{remark}[thm]{Remark}
\newtheorem{example}[thm]{Example}
\numberwithin{equation}{section}

\def\ee{\varepsilon}
\def\qed{{\hfill $\Box$ \bigskip}}
\def\HH{{\mathcal H}}
\def\MM{{\mathcal M}}
\def\BB{{\mathcal B}}
\def\LL{{\mathcal L}}
\def\FF{{\mathcal F}}
\def\EE{{\mathcal E}}
\def\QQ{{\mathcal Q}}
\def\SS{{\mathcal S}}
\def\Cap{{\rm{Cap}}}
\def\wt{\widetilde}
\def\R{{\mathbb R}}
\def\L{{\bf L}}
\def\E{{\mathbb E}}
\def\F{{\bf F}}
\def\P{{\mathbb P}}
\def\H{{\mathbb H}}
\def\N{{\mathbb N}}
\def\eps{\varepsilon}
\def\wh{\widehat}
\def\pf{\noindent{\bf Proof.} }

\title{\Large \bf
Potential Theory of Subordinate Brownian Motions Revisited
}

\author{{\bf Panki Kim}\thanks{This work was supported by Basic
Science Research Program through the National Research Foundation of
Korea(NRF) grant funded by the Korea
government(MEST)(2010-0001984).} \quad {\bf Renming Song} \quad and
\quad {\bf Zoran Vondra{\v{c}}ek}\thanks{Supported in part by the
MZOS grant 037-0372790-2801.} }

\date{ }
\maketitle

\begin{abstract}
The paper discusses and surveys some aspects of the potential theory of subordinate Brownian motion under the assumption that the Laplace exponent of the corresponding subordinator is comparable to a regularly varying function at infinity. This extends some results previously obtained under stronger conditions.
\end{abstract}

\noindent {\bf AMS
2010 Mathematics Subject Classification}: Primary 60J45
Secondary 60G51, 60J35, 60J75

\noindent {\bf Keywords and phrases}: Subordinator, subordinate
Brownian motion, potential theory, Green function, L\'evy density,
Harnack inequality, boundary Harnack principle

\section{Introduction}\label{ksv-sec-intro}

An ${\mathbb R}^d$-valued process $X=(X_t:\, t\ge 0)$ is called a L\'evy
process in ${\mathbb R}^d$ if it is a right continuous process with left
limits and if, for every $s, t\ge 0$, $X_{t+s}-X_s$ is independent
of $\{X_r, r\in [0, s]\}$ and has the same distribution as
$X_s-X_0$. A L\'evy process is completely characterized by its
L\'evy exponent $\Phi$ via
$$
{\mathbb E}[\exp\{i\langle \xi, X_t-X_0 \rangle\}]=\exp\{-t\Phi(\xi)\}, \quad
t\ge 0, \xi\in {\mathbb R}^d.
$$
The L\'evy exponent $\Phi$ of a L\'evy process is given by the
L\'evy-Khintchine formula
$$
\Phi(\xi)=i\langle
l, \xi \rangle + \frac12 \langle \xi, A\xi^T
\rangle + \int_{{\mathbb R}^d}\left(1-e^{i\langle \xi, x \rangle}+i\langle
\xi, x \rangle{\bf 1}_{\{|x|<1\}}\right)\Pi(dx), \quad \xi\in {\mathbb R}^d\, ,
$$
where $ l\in {\mathbb R}^d$, $A$ is a nonnegative definite $d\times d$
matrix, and $\Pi$ is a measure on ${\mathbb R}^d\setminus\{0\}$ satisfying
$\int (1\wedge |x|^2)\, \Pi(dx) <\infty$. $A$ is called the
diffusion matrix, $\Pi$ the L\'evy measure, and
$(l, A, \Pi)$ the generating triplet of the process.

Nowadays L\'evy processes are widely used in various fields, such as
mathematical finance, actuarial mathematics
and  mathematical physics.
However, general L\'evy processes are not very easy to deal with.

A subordinate Brownian motion in ${\mathbb R}^d$ is a L\'evy process which
can be obtained by replacing the time of Brownian motion in ${\mathbb R}^d$ by
an independent subordinator (i.e., an increasing L\'evy process starting from 0).
More precisely, let $B=(B_t:\, t\ge 0)$ be a Brownian motion in
${\mathbb R}^d$ and $S=(S_t:\, t\ge 0)$ be a subordinator independent of $B$.
The process $X=(X_t:\, t\ge 0)$ defined by $X_t=B_{S_t}$ is a
rotationally invariant L\'evy process in ${\mathbb R}^d$ and is called a
subordinate Brownian motion.

The subordinator $S$ used to define the subordinate Brownian motion
$X$ can be interpreted as ``operational'' time or ``intrinsic''
time. For this reason, subordinate Brownian motions have been used
in mathematical finance and other applied fields. Subordinate
Brownian motions form a very large class of L\'evy processes.
Nonetheless, compared with general L\'evy processes, subordinate
Brownian motions are much more tractable. If we take the Brownian
motion $B$ as given, then $X$ is completely determined by the
subordinator $S$. Hence, one can deduce properties of $X$ from
properties of the subordinator $S$. On the analytic level this
translates to the following: Let $\phi$ denote the Laplace exponent
of the subordinator $S$, that is, ${\mathbb E}[\exp\{-\lambda S_t\}]
=\exp\{-t\phi(\lambda)\}$, $\lambda >0$. Then the characteristic
exponent $\Phi$ of the subordinate Brownian motion $X$ takes on the
very simple form $\Phi(x)=\phi(|x|^2)$ (our Brownian motion $B$ runs
at twice the usual speed). Therefore, properties of $X$ should follow
from properties of the Laplace exponent $\phi$.

The Laplace exponent $\phi$ of a subordinator $S$ is a Bernstein
function, hence it has a representation of the form
$$
\phi(\lambda)=b\lambda +\int_{(0,\infty)}(1-e^{-\lambda t})\, \mu(dt)
$$
where $b\ge 0$ and $\mu$ is a measure on $(0,\infty)$ satisfying
$\int_{(0,\infty)}(1\wedge t)\, \mu(dt)<\infty$. If $\mu$ has a
completely monotone density, the function $\phi$ is called a
complete Bernstein function. The purpose of this work is to study
the potential theory of subordinate Brownian motion under the
assumption that the Laplace exponent $\phi$ of the subordinator is a
complete Bernstein function comparable to a regularly varying
functions at infinity.
 More precisely, we will assume that there exist $\alpha\in (0,2)$
and a function $\ell$ slowly varying at infinity such that
$$
\phi(\lambda)\asymp \lambda^{\alpha/2}\ell(\lambda)\, , \quad
\lambda \to \infty\, .
$$
Here and later, for two functions $f$ and $g$ we write
$f(\lambda)\asymp g(\lambda)$ as $\lambda \to \infty$ if the
quotient $f(\lambda)/g(\lambda)$ stays bounded between two positive
constants as $\lambda \to \infty$.

A lot of progress has been made in recent years in the study of the
potential theory of subordinate Brownian motions, see, for instance
\cite{CKSV, CKSV2, KSV1, KSV2, KSV3, RSV, SiSV} and \cite[Chapter
5]{BBKRSV}. In particular, an extensive survey of results obtained
before 2007 is given in \cite[Chapter 5]{BBKRSV}. At that
time, the focus was on the potential theory of the process $X$ in
the whole ${\mathbb R}^d$, the results for (killed) subordinate
Brownian motions in an open subset still being out of reach. In the
last few years significant progress has been made in studying the
potential theory of subordinate Brownian motions killed upon exiting
an open subset of ${\mathbb R}^d$. The main results include the boundary
Harnack principle and sharp Green function estimates. For
processes having a continuous component see \cite{KSV2} (for the
one-dimensional case) and \cite{CKS5, CKSV, CKSV2} (for
multi-dimensional case). For purely discontinuous processes, the
boundary Harnack principle was obtained in \cite{KSV1} and sharp
Green function estimates were discussed in the recent preprint
\cite{KSV3}. The main assumption in \cite{CKSV, CKSV2, KSV1} and
\cite[Chapter 5]{BBKRSV} is that the Laplace exponent of the
subordinator is regularly varying at infinity. The results were
established under different assumptions, some of which turned out to
be too strong and some even redundant. Time is now ripe to put some
of the recent progress under one unified setup and to give a survey
of some of these results. The survey builds upon the work done in
\cite[Chapter 5]{BBKRSV} and \cite{KSV1}. The setup we are going to
assume is more general than all these of the previous papers, so in
this sense, most of the results contained in this paper are
extensions of the existing ones.

In Section \ref{ksv-sec-subordinators} we first recall some basic facts
about subordinators, Bernstein functions and complete Bernstein
functions. Then we establish asymptotic behaviors, near the origin,
of the potential density and L\'evy density of subordinators.

In Section \ref{ksv-sec-sbm} we establish the asymptotic behaviors, near
the origin, of the Green function and the L\'evy density of
our subordinate Brownian motion. These results follow from the
asymptotic behaviors, near the origin, of the potential density and
L\'evy density of the subordinator.

In Section \ref{ksv-sec-hibhp} we prove that the Harnack inequality and
the boundary Harnack principle hold for our subordinate Brownian
motions.

The materials covered in this paper by no means include all that can
be said about the potential theory of subordinate Brownian motions.
One of the omissions is the sharp Green function estimates of
(killed)
subordinate Brownian motions in bounded $C^{1, 1}$ open sets
obtained in the recent preprint \cite{KSV3}. The present paper
builds up the framework for \cite{KSV3} and can be regarded as a
preparation for \cite{KSV3} in this sense.
Another omission is the Dirichlet heat kernel estimates of
subordinate Brownian motions in smooth open sets recently
established in \cite{CKS1, CKS2, CKS3, CKS4}. One of the reasons we
do not include these recent results in this paper is that all these
heat kernel estimates are for particular subordinate Brownian
motions only and are not yet established in the general case.
A third notable omission is the spectral theory for killed
subordinate Brownian motions developed in \cite{CS05, CS06a, CS06c}.
Some of these results have been summarized in \cite[Section
12.3]{SSV}.
A fourth notable omission is the potential theory of subordinate
killed Brownian motions developed in \cite{GPRSSV, GRSS, Son04,
SV03, SV04b, SV06}. Some of these results have been summarized in
\cite[Section 5.5]{BBKRSV} and \cite[Chapter 13]{SSV}. In this paper
we concentrate on subordinate Brownian motions without diffusion
components and therefore this paper does not include results
from\cite{CKSV, CKSV2, KSV2}. One of the reasons for this is that
subordinate Brownian motions with diffusion components require a
different treatment.

We end this introduction with few words on the notations. For
functions $f$ and $g$ we write $f(t)\sim g(t)$ as $t\to 0+$ (resp.
$t\to \infty$) if the quotient $f(t)/g(t)$ converges to 1 as $t\to
0+$ (resp. $t\to \infty$), and $f(t) \asymp g(t)$ as $t\to 0+$
(resp. $t\to \infty$) if the quotient $f(t)/g(t)$ stays bounded
between two positive constants as $t\to 0+$ (resp. $t\to \infty$).

\section{Subordinators}\label{ksv-sec-subordinators}\label{ksv-sec-sub}

\subsection{Subordinators and Bernstein functions}\label{ksv-ss:bf}

Let $S=(S_t:\, t\ge 0)$ be a subordinator, that is, an increasing
L\'evy process taking values in $[0,\infty)$ with $S_0=0$. A
subordinator $S$ is completely characterized by its Laplace exponent
$\phi$ via
$$
{\mathbb E}[\exp(-\lambda S_t)]=\exp(-t \phi(\lambda))\, ,\quad  \lambda > 0.
$$
The Laplace exponent $\phi$ can be written in the form (cf. \cite[p.
72]{Ber})
$$
\phi(\lambda)=b\lambda +\int_0^{\infty}(1-e^{-\lambda t})\,
\mu(dt)\, .
$$
Here $b \ge 0$, and $\mu$ is a $\sigma$-finite measure on
$(0,\infty)$ satisfying
$$
\int_0^{\infty} (t\wedge 1)\, \mu(dt)< \infty\, .
$$
The constant $b$ is called the drift, and $\mu$ the L\'evy measure
of the subordinator $S$.

A $C^{\infty}$ function $\phi:(0,\infty)\to [0,\infty)$ is called a
Bernstein function if $(-1)^n D^n \phi\le 0$ for every positive
integer $n$. Every Bernstein function has a representation (cf.
\cite[Theorem 3.2]{SSV})
$$
\phi(\lambda)=a+b\lambda +\int_{(0,\infty)}(1-e^{-\lambda t})\,
\mu(dt)
$$
where $a,b\ge 0$ and $\mu$ is a measure on $(0,\infty)$ satisfying
$\int_{(0,\infty)}(1\wedge t)\, \mu(dt)<\infty$. $a$ is called the
killing coefficient, $b$ the drift and $\mu$ the L\'evy measure of
the Bernstein function. Thus a nonnegative function $\phi$ on $(0,
\infty)$ is the Laplace exponent of a subordinator if and only if it
is a Bernstein function with $\phi(0+)=0$.

Sometimes we need to deal with killed subordinators, that is,
subordinators killed at independent exponential times.
Let $e_a$ be an exponential random variable with parameter
$a\ge0$, i.e., ${\mathbb P}( e_a >t)=e^{-at}, \, t>0$.
We allow $a=0$ in which case $e_a=\infty$. Assume that $S$ is
a subordinator with Laplace exponent $\phi$ and
$e_a$ is independent of $S$.
We define a process $\widehat S$ by
 $$
 \widehat S_t= \begin{cases}
 S_t , & t < e_a\\
 \infty & t \ge e_a
 \end{cases}.
 $$
The process $\widehat S$ is  the subordinator $S$ killed at an independent exponential time. We call $\widehat S$
a killed subordinator.
The
corresponding Laplace exponent $\widehat  \phi$ is related  to $\phi$ as
$$
\widehat  \phi(\lambda)=a+\phi(\lambda), \qquad \lambda>0.
$$
In fact,
$$
{\mathbb E}[e^{i \xi \cdot \widehat S_t}]= {\mathbb E}[e^{i \xi \cdot S_t} {\bf 1}_{\{  t <   e_a   \}}]={\mathbb E}[e^{i \xi \cdot S_t}]{\mathbb P}(  t <   e_a  )=e^{-t \phi(\lambda)}
e^{-at}=e^{-t(a+\phi(\lambda))}.$$
A function $\phi: (0, \infty)\to [0, \infty)$ is the Laplace
exponent of a killed subordinator if and only if $\phi$ is a
Bernstein function. For this reason, we use a killed subordinator sometimes.

A Bernstein function $\phi$ is called a complete Bernstein function
if the L\'evy measure $\mu$ has a completely monotone density
$\mu(t)$, i.e., $(-1)^n D^n \mu\ge 0$ for every non-negative integer
$n$. Here and below, by abuse of notation  we will denote the L\'evy
density by $\mu(t)$. Complete Bernstein functions form a large
subclass of Bernstein functions. Most of the familiar Bernstein
functions are complete Bernstein functions. See \cite[Chapter
15]{SSV} for an extensive table of complete Bernstein functions.
Here are some examples of complete Bernstein functions:
\begin{description}
\item{(i)} $\phi(\lambda)=\lambda^{\alpha/2}$, $\alpha\in (0, 2]$;
\item{(ii)} $\phi(\lambda)=(\lambda+m^{2/\alpha})^{\alpha/2}-m$, $\alpha\in (0, 2), m\ge
0$;
\item{(iii)} $\phi(\lambda)=\lambda^{\alpha/2}+ \lambda^{\beta/2}$, $0\le
\beta<\alpha\in (0, 2]$;
\item{(iv)} $\phi(\lambda)=\lambda^{\alpha/2}(\log(1+\lambda))^{\gamma/2}$, $\alpha\in (0, 2),
\gamma\in (0, 2-\alpha)$;
\item{(v)} $\phi(\lambda)=\lambda^{\alpha/2}(\log(1+\lambda))^{-\beta/2}$, $0\le
\beta<\alpha\in (0, 2]$.
\end{description}
An example of a Bernstein function which is not a complete Bernstein
function is $1-e^{-\lambda}$.

It is known (cf. \cite[Proposition 7.1]{SSV}) that $\phi$ is a
complete Bernstein function if and only if the function
$\lambda/\phi(\lambda)$ is a complete Bernstein function. For other
properties of complete Bernstein functions we refer the readers to
\cite{SSV}.

The following result, which will play an important role later, says
that the L\'evy density of a complete Bernstein function cannot
decrease too fast in the following sense.

\begin{lemma}[{\cite[Lemma 2.1]{KSV3}}]\label{ksv-l:H2-valid}
Suppose that $\phi$ is a complete Bernstein function with L\'evy
density $\mu$. Then there exists $C_1>0$ such that $\mu(t)\le C_1
\mu(t+1)$ for every  $t>1$.
\end{lemma}

\pf Since $\mu$ is a completely monotone function, by Bernstein's
theorem (cf. \cite[Theorem 1.4]{SSV}), there exists a measure $m$ on
$[0,\infty)$ such that $\mu(t)=\int_{[0,\infty)}e^{-tx} m(dx).$
Choose $r>0$ such that $\int_{[0, r]}e^{-x}\, m(dx)\ge \int_{(r,
\infty)}e^{-x}\, m(dx).$ Then, for any $t>1$, we have
\begin{eqnarray*}
\int_{[0, r]}e^{-t x}\, m(dx)&\ge&e^{-(t -1)r}\int_{[0, r]}e^{-x}\, m(dx)\\
&\ge &e^{-(t -1)r}\int_{(r, \infty)}e^{-x}\, m(dx)\,\ge \, \int_{(r,
\infty)}e^{-t x}\, m(dx).
\end{eqnarray*}
Therefore, for any $t>1$,
\begin{eqnarray*}
&&\mu(t+1)\ge \int_{[0, r]}e^{-(t+1) x}\, m(dx)\ge e^{-r}\int_{[0,
r]}e^{- t x}\, m(dx) \\
&&\ge \frac12\, e^{-r}\int_{[0, \infty)}e^{-t
x}\, m(dx)=\frac12\, e^{-r}\mu(t).
\end{eqnarray*}
\qed

The potential measure of
the (possibly killed) subordinator
 $S$ is defined by
\begin{equation}\label{ksv-potential measure}
U(A)={\mathbb E} \int_0^{\infty}
{\bf 1}_{\{S_t\in A\}}
\, dt, \quad A\subset [0,
\infty).
\end{equation}
Note that $U(A)$ is the expected time the subordinator $S$ spends in the set $A$.
The Laplace transform of the measure $U$ is given by
\begin{equation}\label{ksv-lt potential measure}
{\mathcal L} U(\lambda)=\int_0^{\infty}e^{-\lambda t}\, dU(t)=
{\mathbb E}\int_0^{\infty} \exp(-\lambda S_t)\, dt
=\frac{1}{\phi(\lambda)}\, .
\end{equation}

We call a subordinator $S$ a complete subordinator if its Laplace
exponent $\phi$ is a complete Bernstein function.  The following
characterization of complete subordinators is due to \cite[Remark
2.2]{SV06} (see also \cite[Corollary 5.3]{BBKRSV}).

\begin{prop}\label{ksv-p:2.2}
Let $S$ be a subordinator with  Laplace exponent $\phi$ and
potential measure $U$. Then $\phi$ is a complete Bernstein function
if and only if
$$
U(dt)=c\delta_0(dt)+u(t)dt
$$
for some $c\ge 0$ and completely monotone function $u$.
\end{prop}

In case the constant $c$ in the proposition above is equal to zero,
we will call $u$ the potential density of the subordinator $S$.

An inspection of the argument, given in \cite[Chapter 5]{BBKRSV} or
\cite{SV06}, leading to the proposition above yields the following
two results (cf. \cite[Corollary 5.4 and Corollary 5.5]{BBKRSV} or
\cite[Corollary 2.3 and Corollary 2.4]{SV06}).

\begin{corollary}\label{ksv-c2.3}
Suppose that $S=(S_t:\, t\ge 0)$ is a subordinator whose Laplace exponent
$$
\phi(\lambda)=b\lambda +\int_0^{\infty} (1-e^{-\lambda t})\, \mu(dt)
$$
is a complete Bernstein function with $b>0$ or $\mu(0,
\infty)=\infty$. Then the potential measure $U$ of $S$ has a
completely monotone density $u$.
\end{corollary}

\pf
By
\cite[Corollary 5.4]{BBKRSV} or \cite[Corollary 2.3]{SV06},
if  the drift of the complete subordinator $S$ is zero or the L\'evy measure $\mu$
has infinite mass, then the constant $c$ in Proposition \ref{ksv-p:2.2} is equal to zero
so the potential measure $U$ of $S$ has a density $u$. The completeness
of the density follows directly from Proposition \ref{ksv-p:2.2}.
\qed

\begin{corollary}\label{ksv-c2.4}
Let $S$ be a complete subordinator with Laplace exponent $\phi(\lambda)=\int_0^{\infty} (1-e^{-\lambda t})\mu(dt)$.
Suppose that  the L\'evy measure $\mu$
has infinite mass.
Then the potential measure of a (killed) subordinator with
Laplace exponent $\psi(\lambda):=\lambda/\phi(\lambda)$ has a
completely monotone density $v$ given by
$$
v(t)=\mu(t, \infty).
$$
\end{corollary}

\pf
Since the drift of $S$ is zero and the L\'evy measure $\mu$
has infinite mass, by
\cite[Corollary 5.5]{BBKRSV} or \cite[Corollary 2.4]{SV06}, we have that
$$
\psi(\lambda)=a
+ \int_0^{\infty} (1-e^{-\lambda t})\, \nu(dt)
$$
where $a=\left(\int_0^{\infty} t\mu(t)dt\right)^{-1}$, the  L\'evy measure $\nu$ of $\psi$
has infinite mass and the potential measure of a possibly killed (i.e., $a>0$) subordinator with
Laplace exponent $\psi$ has a
 density $v$ given by
$
v(t)=\mu(t, \infty).
$
The completeness
of the density follows  from \cite[Corollary 5.3]{BBKRSV}, which works for killed subordinators.
\qed

\subsection{Asymptotic behavior of the potential and L\'evy
densities}\label{ksv-ss:asympumu}

 From now on we will always assume that $S$ is a complete
subordinator without drift  and that the Laplace exponent $\phi$ of
$S$ satisfies $\lim_{\lambda\to\infty}\phi(\lambda)=\infty$ (or
equivalently, the L\'evy measure of $S$ has infinite mass). Under
this assumption,  the potential measure $U$ of $S$ has a completely
monotone density $u$ (cf. Corollary \ref{ksv-c2.3}). The main purpose of
this subsection is to determine the asymptotic behaviors of $u$ and
$\mu$ near the origin. For this purpose, we will need the following
result due to Z\"ahle (cf. \cite[Theorem 7]{Z}).

\begin{prop}\label{ksv-p:zahle}
Suppose that $w$ is a
completely monotone function given by
$$
w(t)=\int^\infty_0e^{-st} f(s)\, ds,
$$
where $f$ is a nonnegative decreasing function. Then
$$
f(s)\le \left(1-e^{-1}\right)^{-1} s^{-1}w(s^{-1}), \quad s>0.
$$
If, furthermore, there exist $\delta\in (0, 1)$ and $a, s_0>0$ such
that
\begin{equation}\label{ksv-e:zahle}
w(\lambda t)\le a \lambda^{-\delta} w(t), \quad \lambda\ge 1, t\ge
1/s_0,
\end{equation}
then there exists $C_2=C_2(w,f,a,s_0, \delta)>0$ such that
$$
f(s)\ge C_2 s^{-1}w(s^{-1}), \quad s\le s_0.
$$
\end{prop}

\pf
Using the assumption that $f$ is a nonnegative decreasing function,
we get that, for any $r>0$, we have
\begin{eqnarray*}
w(t)&=&\frac1t\int^\infty_0e^{-s}f\left(\frac{s}{t}\right)ds\\
&\ge&\frac1t\int^r_0e^{-s}f\left(\frac{s}{t}\right)ds
\,\ge\,\frac1tf\left(\frac{r}{t}\right)\left(1-e^{-r}\right).
\end{eqnarray*}
Thus
$$
f\left(\frac{r}{t}\right)\le\frac{tw(t)}{1-e^{-r}}, \quad t>0, r>0.
$$
In particular, we have
$$
f(s)\le \left(1-e^{-1}\right)^{-1}s^{-1}w(s^{-1}), \quad s>0,
$$
and
\begin{equation}\label{ksv-e:zahle1}
f\left(\frac{s}{t}\right)\le
\left(1-e^{-1}\right)^{-1}\frac{t}{s}w\left(\frac{t}{s}\right),
\quad s>0, t>0.
\end{equation}

On the other hand, for $r\in (0, 1]$, we have
\begin{eqnarray*}
tw(t)&=&\int^r_0e^{-s}f\left(\frac{s}{t}\right)ds+\int_r^\infty
e^{-s}f\left(\frac{s}{t}\right)ds\\
&\le&\int^r_0e^{-s}f\left(\frac{s}{t}\right)ds +
f\left(\frac{r}{t}\right)e^{-r}\\
&\le&\left(1-e^{-1}\right)^{-1}t\int^r_0e^{-s}\frac1s\, w\left(\frac{t}{s}\right)ds
+ f\left(\frac{r}{t}\right)e^{-r},
\end{eqnarray*}
where in the last line we used \eqref{ksv-e:zahle1}.
Now we assume \eqref{ksv-e:zahle}, then
we get that
$$
w\left(\frac{t}{s}\right)\le a s^{\delta}w(t), \quad t\ge 1/s_0,
s<r.
$$
Thus, for $r\in (0, 1]$, we have,
$$
tw(t)\le
a\left(1-e^{-1}\right)^{-1}tw(t)\int^r_0e^{-s}s^{\delta-1}ds +
f\left(\frac{r}{t}\right)e^{-r}.
$$
Choosing $r\in(0, 1]$ small enough so that
$$
a\left(1-e^{-1}\right)^{-1}\int^r_0e^{-s}s^{\delta-1}ds\le \frac12,
$$
we conclude that for this choice of $r$, we have
$$
f\left(\frac{r}{t}\right)\ge c_1 tw(t), \quad t\ge 1/s_0
$$
for some constant $c_1>0$. Since $w$ is decreasing, we have
$$
f(s)\ge c_1\frac{r}{s}w\left( \frac{r}{s}\right) \ge
c_2 s^{-1}w(s^{-1}), \quad s\le rs_0,
$$
where $c_2=c_1r$. From this we immediately get that there exists
$c_3>0$ such that
$$
f(s)\ge c_3 s^{-1}w(s^{-1}), \quad s\le s_0.
$$
 \qed

\begin{corollary}
The potential density $u$ of $S$ satisfies
\begin{equation}\label{ksv-e:u-upper-bound}
u(t)\le C_3  t^{-1}\phi(t^{-1})^{-1}\, ,\quad t>0\, .
\end{equation}
\end{corollary}
\pf
Apply the first part of Proposition \ref{ksv-p:zahle} to the
function
$$
w(t):=\int_0^{\infty}e^{-s t}u(s)\, ds =\frac{1}{\phi(t)}.
$$
\qed

We introduce now the main assumption on our Laplace exponent $\phi$
of the complete subordinator $S$ that we will use throughout the
rest of the paper. Recall that a function $\ell:(0,\infty)\to
(0,\infty)$ is slowly varying at infinity if
$$
\lim_{t\to \infty}\frac{\ell(\lambda t)}{\ell(t)}=1\, ,\quad \textrm{for every }\lambda >0\, .
$$

\medskip

\noindent {\bf Assumption (H):} There exist $\alpha\in (0,2)$ and a
function $\ell:(0,\infty)\to (0,\infty)$ which is measurable,
locally bounded above and below by positive constants,
and slowly varying at infinity such that
\begin{equation}\label{ksv-e:reg-var}
\phi(\lambda) \asymp \lambda^{\alpha/2}\ell(\lambda)\, ,\quad \lambda \to \infty\, .
\end{equation}

\begin{remark}\label{ksv-r-interpretation-H}{\rm
The precise interpretation of \eqref{ksv-e:reg-var} will be as follows:
There exists a positive constant $c>1$ such that
$$
c^{-1}\le \frac{\phi(\lambda)}{\lambda^{\alpha/2}\ell(\lambda)} \le
c \qquad \textrm{for all }\lambda \in [1,\infty)\, .
$$
The choice of the interval $[1,\infty)$ is, of course, arbitrary.
Any interval $[a,\infty)$ would do, but with a different constant.
This follows from the
continuity of $\phi$ and the assumption that $\ell$ is
locally bounded above and below by positive constants.
Moreover, by choosing $a>0$ large enough, we could dispense with the
local boundedness assumption. Indeed, by \cite[Lemma 1.3.2]{BGT},
every slowly varying function at infinity is locally bounded on
$[a,\infty)$ for $a$ large enough.

Although the choice of interval $[1,\infty)$ is arbitrary, it will
have as a consequence the fact that all relations of the type
$f(t)\asymp g(t)$ as $t \to \infty$ (respectively $t\to 0+$)
following from \eqref{ksv-e:reg-var} will be interpreted as
$\tilde{c}^{-1} \le f(t)/g(t) \le \tilde{c}$ for $t\ge 1$
(respectively $0<t\le 1$). }
\end{remark}

The assumption \eqref{ksv-e:reg-var} is a very weak assumption on the
asymptotic behavior of $\phi$ at infinity. All the examples (in (i),
(iii) and (v), we need to take $\alpha<2$) above Lemma \ref{ksv-l:H2-valid}
satisfy this assumption. In fact they satisfy the following stronger
assumption
\begin{equation}\label{ksv-e:reg-var2}
\phi(\lambda)= \lambda^{\alpha/2}\ell(\lambda)\, ,
\end{equation}
where $\ell$ is a function slowly varying at infinity. By inspecting
the table in \cite[Chapter 15]{SSV}, one can come up with a lot more
examples of complete Bernstein functions satisfying this stronger
assumption.
In the next example we construct a complete Bernstein function
satisfying \eqref{ksv-e:reg-var}, but not the stronger
\eqref{ksv-e:reg-var2}.

\begin{example}{\rm
Suppose that $\alpha\in (0, 2)$. Let $F$ be a function on $[0,
\infty)$ defined by $F(x)=0$ on $0 \le x<1$
and
$$
F(x)=2^n\, , \quad 2^{2(n-1)/\alpha} \le x < 2^{2n/\alpha},\  n=1, 2,
\dots.
$$
Then clearly $F$ is non-decreasing and  $x^{\alpha/2} \le F(x) \le 2
x^{\alpha/2}$ for all
$x \ge 1$.
 This implies that for all $t>0$,
$$
\frac{t^{\alpha/2}}{2} \le \liminf_{x\to \infty} \frac{F(tx)}{F(x)} \le
\limsup_{x\to \infty} \frac{F(tx)}{F(x)}  \le 2 t^{\alpha/2}.
$$
If $F$ were regularly varying, the above inequality would imply that
the index was $\alpha/2$, hence the limit of $F(tx)/F(x)$ as
$x\to \infty$ would be equal to $c t^{\alpha/2}$ for some positive
constant $c$. But this does not happen because of the following.
Take $t=2^{2/\alpha}$ and a subsequence $x_n=2^{2n/\alpha}$. Then $t
x_n= 2^{2(n+1)/\alpha}$ and therefore
$$
F(t x_n)/F(x_n)=2^{n+2} / 2^{n+1}=2
$$
which should be equal to $c t^{\alpha/2}=c
(2^{2/\alpha})^{\alpha/2}=2c$, implying $c=1$. On the other hand,
take any $t\in (1, 2^{2/\alpha})$ and the same subsequence
$x_n=2^{2n/\alpha}$. Then $t x_n\in [2^{2n/\alpha},
2^{2(n+1)/\alpha} )$ implying $F(t x_n)=F(x_n)$. Thus the quotient
$F(t x_n)/F(x_n)=1$ which should be equal to $c
t^{\alpha}=t^{\alpha}$ for all $t\in (1,2^{1/a})$. Clearly this is
impossible, so $F$ is not regularly varying. This also shows that
$F(x)$ is not $\sim$ to any $c x^{\alpha/2}$, as $x\to \infty$.

Let $\sigma$ be the measure corresponding to the nondecreasing function $F$ (in the sense that $\sigma(dt)=F(dt)$):
$$
\sigma:=\sum^{\infty}_{n=1} 2^n\delta_{2^{2n/\alpha}}\, .
$$
Since $\int_{(0,\infty)}(1+t)^{-1}\, \sigma(dt)<\infty$, $\sigma$ is a Stieltjes measure. Let
$$
g(\lambda):=\int_{(0,\infty)}\frac{1}{\lambda+t}\, \sigma(dt)=\sum_{n=1}^{\infty} \frac{2^n}{\lambda +2^{2n/\alpha}}
$$
be the corresponding Stieltjes function.
It follows from \cite[Theorem 1.7.4]{BGT} or \cite[Lemma 6.2]{WYY}
that $g$ is not regularly varying at infinity.
Moreover, since $F(x)\asymp x^{\alpha/2}$,
$x \to \infty$,
it follows from \cite[Lemma 6.3]{WYY} that $g(\lambda)\asymp \lambda^{\alpha/2-1}$, $\lambda \to \infty$.
Therefore, the function $f(\lambda):=1/ g(\lambda)$ is a
complete Bernstein function which is not regularly varying at infinity,
but satisfies $f(\lambda)\asymp \lambda^{1-\alpha/2}$, $\lambda \to \infty$.
}
\end{example}

Now we are going to establish the asymptotic behaviors of $u$ and
$\mu$ under the assumption {\bf (H)}.

First we claim that under the assumption \eqref{ksv-e:reg-var}, there exist
$\delta\in (0, 1)$ and $a, s_0>0$ such that
\begin{equation}\label{ksv-e:zahle3}
\phi(\lambda t)\ge a\lambda^{\delta}\phi(t), \quad \lambda\ge 1,
t\ge 1/s_0.
\end{equation}
Indeed, by Potter's theorem (cf. \cite[Theorem 1.5.6]{BGT}), for
$0<\epsilon<\alpha/2$ there exists $t_1$ such that
$$
\frac{\ell( t)}{\ell(\lambda t)}\le 2 \max\left(\left(\frac{
t}{\lambda t}\right)^{\epsilon}, \left(\frac{\lambda
t}{t}\right)^{\epsilon}\right)=2\lambda^{\epsilon}\, ,\quad \lambda \ge 1,
t\ge t_1\, .
$$
Hence,
$$
\phi(\lambda t)\ge  c_2 (\lambda t)^{\alpha/2}\ell(\lambda t) = c_2
t^{\alpha/2}\ell(t) \lambda^{\alpha/2} \frac{\ell(\lambda
t)}{\ell(t)}\ge c_3 \phi(t)\lambda^{\alpha/2-\epsilon}\, ,\quad
\lambda\ge 1, t\ge t_1.
$$
By taking $\delta:=\alpha/2 -\epsilon\in (0,1)$, $a=c_3$,  and
$s_0=1/t_1$ we arrive at \eqref{ksv-e:zahle3}.

\begin{thm}\label{ksv-t:behofu}
Let $S$ be a
complete (possibly killed) subordinator
with Laplace exponent $\phi$ satisfying
{\bf (H)}. Then the potential density $u$ of $S$ satisfies
\begin{equation}\label{ksv-e:behofu}
u(t)\asymp
t^{-1}\phi(t^{-1})^{-1}\asymp\frac{t^{\alpha/2-1}}{\ell(t^{-1})}\, ,
\quad t \to 0+\,.
\end{equation}
\end{thm}

\pf
Put
$$
w(t):=\int_0^{\infty}e^{-s t}u(s)\, ds =\frac{1}{\phi(t)},
$$
then by \eqref{ksv-e:zahle3} we have
$$
w(\lambda t)\le a^{-1}\lambda^{-\delta}w(t), \quad \lambda\ge 1, t\ge
1/s_0.
$$
Applying the second part of Proposition \ref{ksv-p:zahle} we see that
there is a constant $c>0$ such that
$$
u(t)\ge c t^{-1}w(t^{-1}),
$$
for small $t>0$. Combining this inequality with
\eqref{ksv-e:u-upper-bound} we arrive at \eqref{ksv-e:behofu}.
 \qed

\begin{thm}\label{ksv-t:behofmu}
Let $S$ be a complete subordinator with Laplace exponent $\phi$
with zero killing coefficient
satisfying {\bf (H)}.
Then the L\'evy density $\mu$ of $S$ satisfies
\begin{equation}\label{ksv-e:behofmu}
\mu(t)\asymp t^{-1}\phi(t^{-1})\asymp t^{-\alpha/2-1}\ell(t^{-1})\,
, \quad t\to 0+\,.
\end{equation}
\end{thm}

\pf  Since $\phi$ is a complete Bernstein function, the function
$\psi(\lambda):=\lambda/\phi(\lambda)$ is also a complete Bernstein
function and satisfies
$$
\psi(\lambda)\asymp \frac{\lambda^{1-\alpha/2}}{\ell(\lambda)}\,
,\quad \lambda \to \infty,
$$
where $\alpha\in (0,2)$ and $\ell$ are the same as in
\eqref{ksv-e:reg-var}. It follows from Corollary \ref{ksv-c2.4} that the
potential measure of
a killed subordinator
 with Laplace exponent $\psi$ has
a complete monotone density $v$ given by
$$
v(t)=\mu(t, \infty)=\int^\infty_t\mu(s)ds.
$$
Applying Theorem \ref{ksv-t:behofu} to $\psi$ and $v$ we get
\begin{equation}\label{ksv-e:v-psi}
\mu(t,\infty)=v(t)\asymp t^{-1}\psi(t^{-1})^{-1}=\phi(t^{-1})\,
,\quad t \to 0\, .
\end{equation}
By using the elementary inequality $1-e^{-c y}\le c(1-e^{-y})$ valid
for all $c \ge 1$ and all $y>0$, we get that $\phi(c\lambda)\le
c\phi(\lambda)$ for all $c\ge 1$ and all $\lambda >0$. Hence
$\phi(s^{-1})=\phi(2 (2s)^{-1})\le 2\phi((2s)^{-1})$ for all $s>0$.
Therefore, by \eqref{ksv-e:v-psi}, for all $s \in (0, 1/2)$
$$
v(s)\le c_1 \phi(s^{-1})\le 2c_1 \phi((2s)^{-1}) \le c_2 v(2s)
$$
for some constants $c_1, c_2>0$. Since
$$
v(t/2)\ge v(t/2)-v(t)=\int_{t/2}^{t}\mu(s)\, ds \ge (t/2)\mu(t)\, ,
$$
we have for all
$t \in (0,1)$,
$$
\mu(t)\le 2 t^{-1}v(t/2)\le c_2 t^{-1} v(t)\le c_3
t^{-1}\phi(t^{-1})\, ,
$$
for some constant $c_3>0$.

Using \eqref{ksv-e:zahle3}
we get that for every $\lambda \ge 1$
$$
\phi(s^{-1})=\phi(\lambda(\lambda s)^{-1})\ge
a\lambda^{\delta}\phi((\lambda s)^{-1})\, ,\quad s\le
\frac{s_0}{\lambda}\, .
$$
It follows from \eqref{ksv-e:v-psi} that there
exists a constant $c_4>0$ such that
$$
c_4^{-1}\phi(s^{-1})\le v(s) \le c_4 \phi(s^{-1})\, , \quad s< 1\,.
$$
Fix $\lambda:=2^{1/\delta}((c_4^2 a^{-1})\vee 1)^{1/\delta}\ge1$.
Then for $s\le (s_0\wedge 1)/\lambda$,
$$
v(\lambda s)\le c_4\phi((\lambda s)^{-1})\le c_4 a^{-1}
\lambda^{-\delta}\phi(s^{-1})\le c_4^2a^{-1}\lambda^{-\delta}
v(s)
\le \frac12 v(s)
$$
by our choice of $\lambda$. Further,
$$
(\lambda-1)s\mu(s)\ge \int_s^{\lambda s}\mu(t)\, dt=v(s)-v(\lambda
s)\ge v(s)-\frac12 v(s)=\frac12 v(s)\, .
$$
This implies that for all small $t$
$$
\mu(t)\ge \frac{1}{2(\lambda-1)}t^{-1} v(t)= c_5 t^{-1}v(t)\ge c_6
t^{-1}\phi(t^{-1})
$$
for some constants $c_5, c_6>0$. The proof is now complete.
\qed

\section{Subordinate Brownian motion}\label{ksv-sec-sbm}

\subsection{Definitions and technical lemma}\label{ksv-ss:sbm}

Let $B=(B_t,
{\mathbb P}_x)$ be a Brownian motion in ${\mathbb R}^d$ with transition
density $p(t,x,y)=p(t,y-x)$ given by
$$
p(t,x)=(4\pi t)^{-d/2}\exp\left(-\frac{|x|^2}{4t}\right), \quad t>0,\,
x,y\in {\mathbb R}^d \, .
$$
The semigroup $(P_t:\, t\ge 0)$ of $B$ is defined by $P_tf(x)=
{\mathbb E}_x[f(B_t)]=\int_{{\mathbb R}^d}p(t,x,y)f(y)\, dy$, where $f$ is a
nonnegative Borel function on ${\mathbb R}^d$. Recall that if $d\ge 3$, the
Green function $G^{(2)}(x,y)=G^{(2)}(x-y)$, $x,y\in {\mathbb R}^d$, of $B$ is
well defined and is equal to
$$
G^{(2)}(x)=\int_0^{\infty}p(t,x)\, dt =
\frac{\Gamma(d/2-1)}{4\pi^{d/2}}\, |x|^{-d+2}\, .
$$

Let $S=(S_t:\, t\ge 0)$ be a complete subordinator independent of $B$, with
Laplace exponent $\phi(\lambda)$, L\'evy measure $\mu$ and potential
measure $U$. In the rest of the paper, we will always assume that
$S$ is a complete subordinator whose killing coefficient is zero, is dependent of $B$ and satisfies ({\bf H}).
Hence $\lim_{\lambda\to \infty}\phi(\lambda)=\infty$,
and thus $S$ has a completely monotone
potential density $u$.
We define a new process $X=(X_t:\, t\ge 0)$ by $X_t:=B_{S_t}$. Then
$X$ is a L\'evy process with characteristic exponent
$\Phi(x)=\phi(|x|^2)$ (see e.g.\cite[pp.197--198]{Sat}) called a
subordinate Brownian motion. The semigroup $(Q_t:\, t\ge 0)$
of the process $X$ is given by
$$
Q_t f(x)= {\mathbb E}_x[f(X_t)]={\mathbb E}_x[f(
B_{S_t})]=\int_0^{\infty} P_s f(x)\,
{\mathbb P}(S_t\in ds)\, .
$$
 The semigroup $Q_t$ has a
density $q(t,x,y)=q(t,x-y)$ given by $q(t,x)=\int_0^{\infty}p(s,x)\,
{\mathbb P}(S_t\in ds)$.

Recall that, according to the criterion of Chung-Fuchs type (see \cite{PS71}
or \cite[pp. 252--253]{Sat}), $X$ is transient if and only if for some small $r>0$,
$\int_{|x|<r} \frac{1}{\Phi(x)}\, dx <\infty$. Since
$\Phi(x)=\phi(|x|^2)$, it follows that $X$ is transient if and only
if
\begin{equation}\label{ksv-transience}
\int_{0+}\frac{\lambda^{d/2-1}}{\phi(\lambda)}\, d\lambda <\infty\, .
\end{equation}
This is always true if $d\ge 3$, and, depending on the subordinator,
may be true for $d=1$ or $d=2$. In the case $d\le 2$, if there
exists $\gamma\in [0, d/2)$ such that
\begin{equation}\label{ksv-e:ass4trans}
\liminf_{\lambda \to 0}\frac{\phi(\lambda)}{\lambda^{\gamma}}>0,
\end{equation}
then \eqref{ksv-transience} holds.

For $x\in {\mathbb R}^d$
and a Borel subset $A$ of ${\mathbb R}^d$, the potential
measure of $X$ is given by
\begin{eqnarray*}
G(x,A)&=&
{\mathbb E}_x\int_0^{\infty}
{\bf 1}_{\{X_t\in A\}}
dt=
\int_0^{\infty}Q_t{\bf 1}_A(x)\, dt
=\int_0^{\infty}\int_0^{\infty}P_s {\bf 1}_A(x){\mathbb P}(S_t\in ds)\, dt\\
&=&\int_0^{\infty}P_s {\bf 1}_A\, u(s)\,ds =\int_A
\int_0^{\infty}p(s,x,y)\, u(s)\, ds\, dy \, ,
\end{eqnarray*}
where the second line follows from (\ref{ksv-potential measure}).
If $X$ is transient and $A$ is bounded, then
$G(x,A)<\infty$ for every
$x\in {\mathbb R}^d$.
In this case
we denote by $G(x,y)$ the density of the potential
measure $G(x,\cdot)$. Clearly, $G(x,y)=G(y-x)$ where
\begin{equation}\label{ksv-green function}
G(x)=\int_0^{\infty} p(t,x)\, U(dt)=\int_0^{\infty} p(t,x) u(t)\,
dt\, .
\end{equation}

The L\'evy measure $\Pi$ of $X$ is given by
(see e.g.~\cite[pp. 197--198]{Sat})
$$
\Pi(A)=\int_A \int_0^{\infty}p(t,x)\, \mu(dt)\, dx =\int_A J(x)\,
dx\, ,\quad A\subset {\mathbb R}^d\, ,
$$
where
\begin{equation}\label{ksv-jumping function}
J(x):= \int_0^{\infty}p(t,x)\, \mu(dt)=\int_0^{\infty}p(t,x)\mu(t)dt
\end{equation}
is the L\'evy density of $X$. Define the function $j:(0,\infty)\to
(0,\infty)$ by
\begin{equation}\label{ksv-function j measure}
j(r):=  \int_0^{\infty} (4\pi)^{-d/2} t^{-d/2} \exp\left(-\frac{r^2}{4t}\right)\,
\mu(dt)\, , \quad r>0\, ,
\end{equation}
and note that by (\ref{ksv-jumping function}), $J(x)=j(|x|)$, $x\in
{\mathbb R}^d\setminus \{0\}$.

Since $x\mapsto p(t,x)$ is continuous and radially decreasing, we conclude that both $G$ and
$J$ are continuous on ${\mathbb R}^d\setminus \{0\}$ and radially decreasing.

The following technical lemma will play a key role in establishing
the asymptotic behaviors of the Green function $G$ and the
L\'evy density $J$ of the subordinate Brownian motion $X$ in the next
subsection.

\begin{lemma}\label{ksv-key technical}
Suppose that $w:(0,\infty)\to (0,\infty)$ is a decreasing function,
$\ell:(0,\infty)\to (0,\infty)$ a measurable
function which is
locally bounded above and below by positive constants and is
slowly varying at $\infty$, and $\beta\in [0,2]$,
$\beta>1-d/2$. If $d=1$ or $d=2$, we additionally assume that there
exist constants $c>0$ and $\gamma <d/2$ such that
\begin{equation}\label{ksv-asymp v2}
w(t)\le ct^{\gamma-1}\, , \quad \forall \, t \ge 1\, .
\end{equation}
Let
$$
I(x):=\int_0^{\infty}(4\pi t)^{-d/2}e^{-\frac{|x|^2}{4t}}w(t)\, dt\, .
$$
\begin{itemize}
    \item[(a)] If
    \begin{equation}\label{ksv-asymp v}
    w(t)\asymp \frac{1}{t^{\beta}\ell(1/t)}\, , \quad t\to 0\, ,
    \end{equation}
    then
    $$
    I(x)\asymp \frac{1}{|x|^{d+2\beta-2}\, \ell \big(\frac{1}{|x|^2}\big)} \asymp\frac{w(|x|^2)}{|x|^{d-2}}\, , \quad |x|\to 0 \, .
    $$
    \item[(b)] If
    \begin{equation}\label{ksv-asymp v-sim}
    w(t)\sim \frac{1}{t^{\beta}\ell(1/t)}\, , \quad t\to 0\, ,
    \end{equation}
    then
    $$
    I(x)\sim \frac{\Gamma(d/2+\beta-1)}{4^{1-\beta}\pi^{d/2}}\, \frac{1}{|x|^{d+2\beta-2}\ell\big(\frac{1}{|x|^2}\big)}\, , \quad |x|\to 0\, .
    $$
\end{itemize}
\end{lemma}
\pf (a) Let us first note that the assumptions of the lemma guarantee
that $I(x)<\infty $ for every $x\neq 0$. Now, let $\xi\ge 1/4$ to be
chosen later. By a change of variable we get
\begin{eqnarray}
\int_0^{\infty}(4\pi t)^{-d/2}e^{-\frac{|x|^2}{4t}}w(t)\, dt
&=& \frac{1}{4\pi^{d/2}}\left(|x|^{-d+2}\int_0^{\xi |x|^2} t^{d/2-2}
e^{-t} w\left(\frac{|x|^2}{4t}\right)\, dt \right.\nonumber\\
&&\qquad \left. +
|x|^{-d+2}\int_{\xi|x|^2}^{\infty} t^{d/2-2} e^{-t}
w\left(\frac{|x|^2}{4t}\right)\, dt\right)\nonumber\\
&=:&
\frac{1}{4\pi^{d/2}}\left(|x|^{-d+2}I_1(x)+|x|^{-d+2}I_2(x)\right)\, .\label{ksv-e:Ione+Itwo}
\end{eqnarray}
We first consider $I_1(x)$ for the case $d=1$ or $d=2$. It follows
from the assumptions that there exists a positive constant $c_1$
such that $w(s)\le c_1 s^{\gamma-1}$ for all $s\ge 1/(4\xi)$. Thus
\begin{eqnarray*}
I_1(x)\le \int_0^{\xi |x|^2}
t^{d/2-2}e^{-t}c_1\left(\frac{|x|^2}{4t}\right)^{\gamma-1}\, dt \le
c_2 |x|^{2\gamma-2}\int_0^{\xi|x|^2}t^{d/2-\gamma-1}\, dt = c_3
|x|^{d-2}\, .
\end{eqnarray*}
It follows that
\begin{equation}\label{ksv-Ione}
\lim_{|x|\to 0} |x|^{-d+2}I_1(x) \left(|x|^{d-2+2\beta}\,
\ell\left(\frac{1}{|x|^2}\right) \right) =0\, .
\end{equation}
In the case $d\ge 3$, we proceed similarly, using the bound $w(s)\le
w(1/(4\xi))$ for $s\ge 1/(4\xi)$.

Now we consider $I_2(x)$:
\begin{eqnarray*}
&&|x|^{-d+2}I_2(x)\,=\,\frac{1}{|x|^{d-2}}\int_{\xi |x|^2}^{\infty}
t^{d/2-2}e^{-t}w\left(\frac{|x|^2}{4t}\right)\, dt\\
&&=\frac{4^\beta}{|x|^{d+2\beta-2}\, \ell(\frac{1}{|x|^2})}\int_{\xi
|x|^2}^{\infty} t^{d/2-2+\beta}e^{-t} w\left(\frac{|x|^2}{4t}\right)
\left(\frac{|x|^2}{4t}\right)^{\beta}\, \ell\,
\left(\frac{4t}{|x|^2}\right)\, \frac{\ell(\frac{1}{|x|^2})}{\ell
(\frac{4t}{|x|^2})}\, dt\, .
\end{eqnarray*}
Using the assumption (\ref{ksv-asymp v}), we can see that there is a
constant $c_1>1$ such that
$$
c_1^{-1} \le w\left(\frac{|x|^2}{4t}\right)
\left(\frac{|x|^2}{4t}\right)^{\beta}\, \ell\,
\left(\frac{4t}{|x|^2}\right) <c_1\, ,
$$
for all $t$ and $x$ satisfying $|x|^2/(4t)\le 1/(4\xi)$.

Now choose a $\delta\in (0,d/2-1+\beta)$ (note that by assumption,
$d/2-1+\beta>0$). By Potter's theorem (cf. \cite[Theorem 1.5.6
(i)]{BGT}), there exists $\rho=\rho(\delta)\ge1$ such that
\begin{equation}\label{ksv-e:potter1}
\frac{\ell(\frac{1}{|x|^2})}{\ell(\frac{4t}{|x|^2})}\le
2\left(\left(\frac{1/|x|^2}{4t/|x|^2}\right)^{\delta}\vee
\left(\frac{1/|x|^2}{4t/|x|^2}\right)^{-\delta}\right)=2\left((4t)^{\delta}\vee
(4t)^{-\delta}\right)\le c_2(t^{\delta}\vee t^{-\delta})
\end{equation}
whenever $\frac{1}{|x|^2}>\rho$ and $\frac{4t}{|x|^2} >\rho$.
By reversing the roles of $1/|x|^2$ and $4t/|x|^2$ we also get that
\begin{equation}\label{ksv-e:potter2}
\frac{\ell(\frac{1}{|x|^2})}{\ell(\frac{4t}{|x|^2})}\ge c_2^{-1}(t^{\delta}\wedge t^{-\delta})\
\end{equation}
for  $\frac{1}{|x|^2}>\rho$ and $\frac{4t}{|x|^2} >\rho$.
Now we define $\xi:=\frac{\rho}{4}$ so that for all $x\neq 0$
with $|x|^2 \le \frac{1}{4\xi}$ and $t>\xi|x|^2$ we have that
\begin{eqnarray}
c_1^{-1} c_2^{-1} \, t^{d/2-2+\beta}e^{-t}
(t^{\delta}\wedge t^{-\delta}) &\le & t^{d/2-2+\beta}e^{-t} w\left(\frac{|x|^2}{4t}\right)
\left(\frac{|x|^2}{4t}\right)^{\beta}\, \ell\,
\left(\frac{4t}{|x|^2}\right)
\frac{\ell(\frac{1}{|x|^2})}{\ell(\frac{4t}{|x|^2})} \,
\nonumber \\
&  \le & c_1 c_2 \, t^{d/2-2+\beta}e^{-t}
(t^{\delta}\vee t^{-\delta})\, .\label{ksv-e:key-lemma-1}
\end{eqnarray}
Let
\begin{eqnarray*}
c_3&:=&c_1^{-1} c_2^{-1}\int_0^{\infty}t^{d/2-2+\beta}e^{-t}
(t^{\delta}\wedge t^{-\delta})
 dt <\infty\, ,\\
c_4&:=& c_1 c_2 \int_0^{\infty} t^{d/2-2+\beta}e^{-t}
(t^{\delta}\vee t^{-\delta})
 dt <\infty\, .
\end{eqnarray*}
The integrals are finite because of assumption $d/2-2+\beta-\delta>-1$.
It follows from \eqref{ksv-e:key-lemma-1} that
\begin{eqnarray*}
c_3&\le & \liminf_{|x|\to 0}\int_0^{\infty} t^{d/2-2+\beta}e^{-t}
w\left(\frac{|x|^2}{4t}\right)
\left(\frac{|x|^2}{4t}\right)^{\beta}\, \ell\,
\left(\frac{4t}{|x|^2}\right)
\frac{\ell(\frac{1}{|x|^2})}{\ell(\frac{4t}{|x|^2})} \,
{\bf 1}_{(\xi|x|^2,\infty)}(t)
dt\\
&\le & \limsup_{|x|\to 0}\int_0^{\infty} t^{d/2-2+\beta}e^{-t} w\left(\frac{|x|^2}{4t}\right)
\left(\frac{|x|^2}{4t}\right)^{\beta}\, \ell\,
\left(\frac{4t}{|x|^2}\right)
\frac{\ell(\frac{1}{|x|^2})}{\ell(\frac{4t}{|x|^2})} \,
{\bf 1}_{(\xi|x|^2,\infty)}(t)
dt \le c_4\, .
\end{eqnarray*}
This means that
\begin{eqnarray}
\lefteqn{|x|^{-d+2}I_2(x) \left(|x|^{d-2\beta+2}\, \ell(\frac{1}{|x|^2}) \right)}\nonumber \\
&=&4^{\beta}\int_{\xi
|x|^2}^{\infty} t^{d/2-2+\beta}e^{-t} w\left(\frac{|x|^2}{4t}\right)
\left(\frac{|x|^2}{4t}\right)^{\beta}\, \ell\,
\left(\frac{4t}{|x|^2}\right)\, \frac{\ell(\frac{1}{|x|^2})}{\ell
(\frac{4t}{|x|^2})}\, dt \asymp 1\, .\label{ksv-Itwo}
\end{eqnarray}
Combining \eqref{ksv-Ione} and \eqref{ksv-Itwo} we have proved the first part of the lemma.

\noindent
(b) The proof is almost the same with a small difference at the very end. Since $\ell$ is slowly varying at $\infty$, we have that
$$
\lim_{|x|\to
0}\frac{\ell(\frac{1}{|x|^2})}{\ell(\frac{4t}{|x|^2})}=1\, .
$$
This implies that
\begin{eqnarray*}
& &\lim_{|x|\to 0} t^{d/2-2+\beta}e^{-t} w\left(\frac{|x|^2}{4t}\right)
\left(\frac{|x|^2}{4t}\right)^{\beta}\, \ell\,
\left(\frac{4t}{|x|^2}\right)
\frac{\ell(\frac{1}{|x|^2})}{\ell(\frac{4t}{|x|^2})} \,
{\bf 1}_{(\xi|x|^2,\infty)}(t)\\
& &\quad = t^{d/2-2+\beta}e^{-t} {\bf 1}_{(0,\infty)}(t)  \, .
\end{eqnarray*}
By the right-hand side inequality in \eqref{ksv-e:key-lemma-1}, we can apply the dominated convergence theorem to conclude that
\begin{eqnarray*}
\lefteqn{\lim_{|x|\to 0}|x|^{-d+2}I_2(x) \left(|x|^{d-2\beta+2}\, \ell(\frac{1}{|x|^2}) \right)}\\
&=&\lim_{|x|\to 0} 4^{\beta}\int_{0}^{\infty} t^{d/2-2+\beta}e^{-t} w\left(\frac{|x|^2}{4t}\right)
\left(\frac{|x|^2}{4t}\right)^{\beta}\, \ell\,
\left(\frac{4t}{|x|^2}\right)\, \frac{\ell(\frac{1}{|x|^2})}{\ell
(\frac{4t}{|x|^2})}{\bf 1}_{(\xi|x|^2,\infty)}(t)\, dt \\
&=&4^{\beta}\Gamma(d/2-1+\beta)\, .
\end{eqnarray*}
Together with \eqref{ksv-e:Ione+Itwo} and \eqref{ksv-Ione} this proves the second part of the lemma.
\qed

\subsection{Asymptotic behavior of the Green function and L\'evy
density}\label{ksv-ss:asympgj}

The goal of this subsection is to establish the asymptotic behaviors
of the Green function $G(x)$ and L\'evy density $J(x)$ of the
subordinate process
$X$ under certain assumptions on the Laplace
exponent $\phi$ of the subordinator $S$. We start with the Green
function.

\begin{thm}\label{ksv-t:Gorigin}
Suppose that the Laplace exponent $\phi$ is a complete Bernstein function satisfying the assumption
{\bf (H)} and that $\alpha\in (0,2\wedge d)$. In the case $d\le 2$,
we further assume \eqref{ksv-e:ass4trans}. Then
$$
G(x)\asymp \frac1{|x|^{d}\phi(|x|^{-2})}\asymp
\frac1{|x|^{d-\alpha}\ell(|x|^{-2})}, \qquad |x|\to 0.
$$
\end{thm}

\pf It follows from Theorem
\ref{ksv-t:behofu}
that the potential
density $u$ of $S$ satisfies
$$
u(t)\asymp
t^{-1}\phi(t^{-1})^{-1}\asymp\frac{t^{\alpha/2-1}}{\ell(t^{-1})}\, ,
\quad t \to 0+\,.
$$
Using \eqref{ksv-e:u-upper-bound} and \eqref{ksv-e:ass4trans}
we conclude that in case $d\le 2$  there exists $c>0$ such that
$$
u(t)\le ct^{\gamma-1}, \quad t\ge 1\, .
$$
We can now apply Lemma \ref{ksv-key technical} with $w(t)=u(t)$,
$\beta=1-\alpha/2$ to obtain the required asymptotic behavior. \qed

\begin{remark}\label{ksv-cond&zahle}{\rm (i) Since $\alpha$ is always assumed to be in $(0, 2)$,
the assumption $\alpha\in (0, 2\wedge d)$ in the theorem above makes a difference only
in the case $d=1$.

\noindent
(ii) In case $d\ge 3$, the conclusion of the theorem above is proved in \cite[Theorem 1 (ii)--(iii)]{Z}
under weaker assumptions. The statement of \cite[Theorem 1 (ii)]{Z} in case $d\le 2$ is incorrect and the proof has an error.
}
\end{remark}

The asymptotic behavior near the origin of $J(x)$ is contained in
the following result.

\begin{thm}\label{ksv-t:Jorigin}
Suppose that the Laplace exponent $\phi$ is a complete Bernstein
function satisfying the assumption {\bf (H)}.
Then
$$
J(x)\asymp \frac{\phi(|x|^{-2})}{|x|^{d}}\asymp
\frac{\ell(|x|^{-2})}{|x|^{d+\alpha}}, \qquad |x|\to 0.
$$
\end{thm}

\pf It follows from Theorem \ref{ksv-t:behofmu} that the
L\'evy density $\mu$ of $S$ satisfies
$$
\mu(t)\asymp t^{-1}\phi(t^{-1})\asymp t^{-\alpha/2-1}\ell(t^{-1})\,
, \quad t \to 0+\,.
$$
Since $\mu(t)$ is decreasing and integrable at infinity,
one can easily show that there exists $c>0$ such that
$$
\mu(t)\le ct^{-1}, \quad t\ge 1.
$$
In fact, if the claim above were not valid, we could find an increasing sequence
$\{t_n\}$ such that $t_1>1, t_n\uparrow\infty, t_n-t_{n-1}\ge t_n/2$ and that
$\mu(t_n)\ge n t_n^{-1}$. Then we would have
$$
\int^\infty_1\mu(t)dt=\int^{t_1}_1\mu(t)dt+\sum^\infty_{n=2}\int^{t_n}_{t_{n-1}}\mu(t)dt
\ge \frac{t_1-1}{t_1}+\sum^\infty_{n=2}\frac{n}2=\infty,
$$
contradicting the integrability of $\mu$ at infinity. Therefore the claim above is valid.
We can now apply Lemma \ref{ksv-key technical} with $w(t)=\mu(t)$,
$\beta=1+\alpha/2$ and $\gamma=0$ to obtain the required asymptotic
behavior. \qed

\begin{prop}\label{ksv-properties of j}
Suppose that the Laplace exponent $\phi$ is a complete Bernstein function satisfying the assumption
{\bf (H)}. Then the following assertions hold.
\begin{description}
\item{(a)} For any $K>0$, there exists $C_4=C_4(K)
>1$ such that
\begin{equation}\label{ksv-H:1}
j(r)\le C_4\, j(2r), \qquad \forall r\in (0, K).
\end{equation}
\item{(b)} There exists $C_5
>1$ such that
\begin{equation}\label{ksv-H:2}
j(r)\le C_5\, j(r+1), \qquad \forall r>1.
\end{equation}
\end{description}
\end{prop}

\pf
\eqref{ksv-H:1} follows immediately from Theorem \ref{ksv-t:Jorigin}.
However, we give below a proof of both \eqref{ksv-H:1} and \eqref{ksv-H:2}
using only \eqref{ksv-nu 0}--\eqref{ksv-nu infty}.

For simplicity we redefine in this proof the function $j$ by
dropping the factor $(4\pi)^{-d/2}$ from its definition. This does
not effect \eqref{ksv-H:1} and \eqref{ksv-H:2}. It follows from Lemma
\ref{ksv-l:H2-valid} and Theorem \ref{ksv-t:behofmu} that
\begin{description}
\item{(a)} For any $K>0$, there exists $c_1=c_1(K)
>1$ such that
\begin{equation}\label{ksv-nu 0}
\mu(r)\le c_1\, \mu(2r), \qquad \forall r\in (0, K).
\end{equation}
\item{(b)} There exists $c_2
>1$ such that
\begin{equation}\label{ksv-nu infty}
\mu(r)\le c_2\, \mu(r+1), \qquad \forall r>1.
\end{equation}
\end{description}
Let $0< r < K$. We have
\begin{eqnarray*}
\lefteqn{j(2r) = \int_0^{\infty} t^{-d/2} \exp(- r^2/t)\mu(t)\, dt} \\
& \ge & \frac12 \left(\int_{K/2}^{\infty}t^{-d/2} \exp(-
r^2/t)\mu(t)\, dt
+\int_0^{2K} t^{-d/2} \exp(- r^2/t)\mu(t)\, dt\right)\\
&=& \frac12 (I_1 + I_2).
\end{eqnarray*}
Now,
\begin{eqnarray*}
I_1 &=& \int_{K/2}^{\infty}t^{-d/2} \exp(- \frac{r^2}{t})\mu(t)\, dt
 = \int_{K/2}^{\infty}t^{-d/2} \exp(-\frac{ r^2}{4t})
\exp(-\frac{3 r^2}{4t})\mu(t)\, dt \\
 &\ge &\int_{K/2}^{\infty}t^{-d/2} \exp(-\frac{ r^2}{4t})
\exp(-\frac{3 r^2}{2K})\mu(t)\, dt
\ge  e^{-3K/2}\int_{K/2}^{\infty}t^{-d/2} \exp(-\frac{ r^2}{4t})
\mu(t)\, dt\, ,\\
& & \\
I_2 &=& \int_0^{2K} t^{-d/2}\exp(- \frac{r^2}{t})\mu(t)\, dt
 = 4^{-d/2+1} \int_0^{K/2} s^{-d/2}\exp(-\frac{ r^2}{4s})\mu(4s)\, ds \\
&\ge & c_1^{-2} 4^{-d/2+1}\int_0^{K/2} s^{-d/2}
\exp(-\frac{ r^2}{4s})\mu(s)\,ds. \\
\end{eqnarray*}
Combining the three displays above we get that $j(2r)\ge c_3\, j(r)$
for all $ r\in (0, K)$.

To prove \eqref{ksv-H:2} we first note that for all $t\ge 2$ and all
$r\ge 1$ it holds that
$$
\frac{(r+1)^2}{t}-\frac{r^2}{t-1}\le 1\, .
$$
This implies that
\begin{equation}\label{ksv-a}
\exp\left(-\frac{(r+1)^2}{4t}\right) \ge e^{-1/4} \exp\left(-\frac{ r^2}{4(t-1)}\right), \quad
\mbox{ for all }r>1, t>2\, .
\end{equation}
Now we have
\begin{eqnarray*}
\lefteqn{j(r+1)=  \int_0^{\infty} t^{-d/2}\exp(-\frac{(r+1)^2}{4t}) \mu(t)\, dt }\\
&\ge & \frac12 \left( \int_0^8 t^{-d/2}\exp(-\frac{(r+1)^2}{4t})
\mu(t)\, dt +\int_3^{\infty}t^{-d/2}\exp(-\frac{(r+1)^2}{4t})
\mu(t)\, dt \right)\\
 &=&\frac12(I_3+I_4).
\end{eqnarray*}
For $I_3$ note that $(r+1)^2\le 4r^2$ for all $r>1$. Thus
\begin{eqnarray*}
I_3 &=& \int_0^8 t^{-d/2}\exp(-\frac{(r+1)^2}{4t}) \mu(t)\, dt
 \ge  \int_0^8 t^{-d/2}\exp(- r^2/t) \mu(t)\, dt \\
&=& 4^{-d/2+1}\int_0^2 s^{-d/2}\exp(-\frac{ r^2}{4s}) \mu(4s)\, ds \\
& \ge &  c_1^{-2} 4^{-d/2+1}\int_0^2 s^{-d/2}\exp(-\frac{ r^2}{4s}) \mu(s)\, ds\, ,\\
& &\\
I_4&=& \int_3^{\infty}t^{-d/2}\exp(-\frac{(r+1)^2}{4t}) \mu(t)\, dt\\
& \ge &  \int_3^{\infty}t^{-d/2}\exp\{-1/4\} \exp(-\frac{ r^2}{4(t-1)})\, \mu(t)\, dt\\
&=&e^{-1/4} \int_2^{\infty} (s-1)^{-d/2}\exp(-\frac{ r^2}{4s})
\, \mu(s+1)\, ds\\
&\ge &  c_1^{-1} e^{-1/4} \int_2^{\infty} s^{-d/2}\exp(-\frac{ r^2}{4s})
\mu(s)\, ds\, .
\end{eqnarray*}
Combining the three displays above we get that
$j(r+1)\ge c_4\, j(r)$ for all $ r>1$.
\qed

\subsection{Some results on subordinate Brownian motion in
${\mathbb R}$}\label{ksv-ss:1dsbm}

In this
subsection
 we assume $d=1$. We will consider subordinate
Brownian motions in ${\mathbb R}$. Let $B=(B_t:\, t\ge 0)$ be a Brownian
motion in ${\mathbb R}$, independent of $S$, with
$$
{\mathbb E}\left[e^{i\theta(B_t-B_0)}\right] =e^{-t\theta^2}, \qquad\,
\forall \theta\in {\mathbb R}, t>0.
$$
The subordinate Brownian motion $X=(X_t:t\ge 0)$ in ${\mathbb R}$ defined by
$X_t=B_{S_t}$ is a symmetric L\'evy process with the characteristic
exponent $\Phi(\theta)=\phi(\theta^2)$ for all $\theta\in {\mathbb R}.$
In the first part of this
subsection,
up to Corollary \ref{ksv-c:vsm}, we do
not need to assume that $\phi$ satisfies the assumption ({\bf H}).

Let $\overline{X}_t:=\sup\{0\vee X_s:0\le s\le t\}$ and let
$L=(L_t:\, t\ge 0)$ be a local time of $\overline{X}-X$ at $0$. $L$
is also called a local time of the process $X$ reflected at the
supremum. Then the right continuous inverse $L^{-1}_t$ of $L$ is a
subordinator and is called the ladder time process of $X$. The
process $\overline{X}_{L^{-1}_t}$ is also a subordinator and is
called the ladder height process of $X$. (For the basic properties
of the ladder time and ladder height processes, we refer our readers
to \cite[Chapter 6]{Ber}.) Let $\chi$ be the Laplace exponent of the
ladder height process of $X$. It follows from \cite[Corollary
9.7]{Fris} that
\begin{equation}\label{ksv-e:formula4leoflh}
\chi(\lambda)=
\exp\left(\frac1\pi\int^{\infty}_0\frac{\log(\Phi(\lambda\theta))}
{1+\theta^2}d\theta \right)
=\exp\left(\frac1\pi\int^{\infty}_0\frac{\log(
 \phi(\lambda^2\theta^2))}{1+\theta^2}d\theta
\right), \quad \forall \lambda>0.
\end{equation}
The next result, first proved independently in \cite{KSV3} and
\cite{Kw}, tells us that $\chi$ is also complete Bernstein function.
The proof presented below is taken from \cite{KSV3}.

\begin{prop}\label{ksv-p:chi is cbf}
Suppose $\phi$, the Laplace exponent of the subordinator $S$, is a
complete Bernstein function. Then the Laplace exponent $\chi$ of the
ladder height process of the subordinate Brownian motion
$X_t=B_{S_t}$ is also  a complete Bernstein function.
\end{prop}
\pf
It follows from Theorem \cite[Theorem 6.10]{SSV} that $\phi$ has
the following representation:
\begin{equation}\label{ksv-e:exp-repr}
\log \phi(\lambda)=\gamma
+\int^\infty_0\left(\frac{t}{1+t^2}-\frac1{\lambda+t}\right)\eta(t)dt,
\end{equation}
where $\eta$ in a function such that $0\le \eta(t)\le 1$ for all
$t>0$.
By \eqref{ksv-e:exp-repr} and \eqref{ksv-e:formula4leoflh},  we have
\begin{eqnarray*}
\log \chi(\lambda)
=\frac{\gamma}2+\frac1{\pi}\int^\infty_0\int^\infty_0
\left(\frac{t}{1+t^2}-\frac1{\lambda^2\theta^2+t}\right)\eta(t)dt
\frac{d\theta}{1+\theta^2}.
\end{eqnarray*}
By using $0\le \eta(t)\le 1$, we have
\begin{eqnarray*}
\eta(t)\left|\frac{t}{1+t^2}-\frac1{\lambda^2\theta^2+t}\right|
\frac{1}{1+\theta^2}&\le& \frac{1}{1+t^2} \frac{1}{1+\theta^2}\left(\frac1{\lambda^2
\theta^2+t }+ \frac{\lambda^2 \theta^2t}{\lambda^2 \theta^2+t
}\right)\\
&  \le& \frac{1}{1+t^2} \left(\frac1{\lambda^2
\theta^2+t }+ \frac{\lambda^2 t}{\lambda^2 \theta^2+t
}\right).
\end{eqnarray*}
Since
\begin{eqnarray*}
\int^\infty_0\frac1{\lambda^2 \theta^2+t }\, d\theta
=\frac1t\int^\infty_0\frac1{\frac{\lambda^2\theta^2}t+1}\, d\theta=\frac1t\frac{\sqrt{t}}{\lambda}\int^\infty_0\frac1{\gamma^2+1}\, d\gamma
= \frac{\pi}{2\lambda\sqrt{t}},
\end{eqnarray*}
we can use  Fubini's theorem to get
\begin{eqnarray}
\log \chi(\lambda)
&=&\frac{\gamma}2+\int^\infty_0\left(\frac{t}{2(1+t^2)}-
\frac1{2\sqrt{t}(\lambda+\sqrt{t})}\right)\eta(t)dt\label{ksv-e:logexp4chi}\\
&=&\frac{\gamma}2+\int^\infty_0\left(\frac{t}{2(1+t^2)}-
\frac1{2(1+t)}\right)\eta(t)dt
\nonumber\\
& & +\int^\infty_0
\left(\frac1{2(1+t)}-\frac1{2\sqrt{t}(\lambda+\sqrt{t})}\right)\eta(t)dt\nonumber\\
&=&\gamma_1+\int^\infty_0\left(\frac{s}{1+s^2}-\frac1{\lambda+s}
\right)\eta(s^2)ds.\nonumber
\end{eqnarray}
Applying \cite[Theorem 6.10]{SSV} we get that $\chi$ is a complete
Bernstein function. \qed

The potential measure of the ladder height process of $X$ is denoted
by $V$ and its density by $v$. We will also use $V$ to denote the
renewal function of $X$: $V(t):=V((0,t))=\int_0^t v(s)\, ds$.

The following result is first proved in \cite{KSV3}.

\begin{prop}\label{ksv-p:chiphi}
$\chi$ is related to $\phi$ by the following relation
$$
e^{-\pi/2}\sqrt{\phi(\lambda^2)}\le \chi(\lambda)\le e^{\pi/2}\sqrt{\phi(\lambda^2)}\, ,\qquad
\textrm{for all }\lambda>0\, .$$
\end{prop}

\pf
According to \eqref{ksv-e:logexp4chi}, we have
$$
\log \chi(\lambda)
=\frac{\gamma}2+\frac12\int^\infty_0\left(\frac{t}{1+t^2}-
\frac1{\sqrt{t}(\lambda+\sqrt{t})}\right)\eta(t)dt\, .
$$
Together with representation \eqref{ksv-e:exp-repr} we get that for all $\lambda >0$
\begin{eqnarray*}
\lefteqn{\left|\, \log \chi(\lambda)-\frac12 \log \phi(\lambda^2)\, \right| }\\
&=&
\frac12\left|\int_0^{\infty}\left(\Big(\frac{t}{1+t^2}- \frac1{\sqrt{t}
(\lambda+\sqrt{t})}\Big) -\Big(\frac{t}{1+t^2}-\frac1{\lambda^2+t}
\Big)\right)\eta(t)\, dt\right|\\
&\le &\frac12\int_0^{\infty}\frac{\lambda(\sqrt{t}+\lambda)}{(\lambda^2+t)
\sqrt{t}(\lambda+\sqrt{t})}\, dt=\frac12\int_0^{\infty}\frac{\lambda}{(\lambda^2+t)\sqrt{t}}\, dt =\frac{\pi}{2}\,
.
\end{eqnarray*}
This implies that
$$
-\pi/2 \le \log \chi(\lambda)-\frac12 \log \phi(\lambda^2) \le \pi/2\,
,\qquad \textrm{for all }\lambda>0\, ,
$$
i.e.,
$$
e^{-\pi/2}\le \chi(\lambda)\phi(\lambda^2)^{-1/2}\le e^{\pi/2}\, ,\qquad
\textrm{for all }\lambda>0\, .
$$
\qed

Combining the above two propositions with Corollary \ref{ksv-c2.3}, we obtain

\begin{corollary}\label{ksv-c:vsm}
Suppose $\phi$, the Laplace exponent of the subordinator $S$, is a
complete Bernstein function satisfying $\lim_{\lambda\to\infty}\phi(\lambda)=\infty$.
Then the potential measure  of the
ladder height process of the subordinate Brownian motion
$X_t=B_{S_t}$ has a completely monotone density $v$. In particular,
$v$ and the renewal function $V$ are $C^{\infty}$ functions.
\end{corollary}

In the remainder of  this paper we will always assume that $\phi$ satisfies
the assumption {\bf(H)}.
We will not explicitly mention this assumption anymore.

\medskip

 Since $\phi(\lambda)\asymp\lambda^{\alpha/2}\ell(\lambda)$ as $\lambda \to
\infty$, Lemma \ref{ksv-p:chiphi} implies that
\begin{equation}\label{ksv-e:abofkappaatinfty}
\chi(\lambda)\asymp\lambda^{\alpha/2}(\ell(\lambda^2))^{1/2}, \qquad
t\to \infty.
\end{equation}
It follows from \eqref{ksv-e:abofkappaatinfty} that $\lim_{\lambda\to
\infty} \chi(\lambda)/\lambda=0$, hence the ladder height process
does not have a drift. Recall that $V(t)=V((0,t))=\int_0^t v(s)ds$
is the renewal function of the ladder height process of $X$. In
light of \eqref{ksv-e:abofkappaatinfty}, we have, as a consequence of
Theorem \ref{ksv-t:behofu}, the following result.

\begin{prop}\label{ksv-p:abofgf4lhpat0}
As $t\to 0$, we have
$$
V(t)\,\asymp
\phi(t^{-2})^{-1/2}\asymp\,\frac{t^{\alpha/2}}{(\ell(t^{-2}))^{1/2}}
$$
and
$$
v(t)\,\asymp
t^{-1}\phi(t^{-2})^{-1/2}\asymp\,\frac{t^{\alpha/2-1}}{ (\ell(t^{-2}))^{1/2}}\, .
$$
\end{prop}

\begin{remark}\label{ksv-r:abofgf4lhpat0}{\rm
It follows immediately from the proposition above that there exists
a positive constant $c>0$ such that $V(2t)\le c V(t)$ for all $t\in
(0,2)$. }
\end{remark}

It follows from \eqref{ksv-e:abofkappaatinfty} above and
\cite[Lemma 7.10]{Ky} that the process $X$ does not creep upwards.
Since $X$ is symmetric, we know that $X$ also does not creep
downwards. Thus if, for any $a\in {\mathbb R}$, we define
$$
\tau_a=\inf\{t>0: X_t<a\}, \quad \sigma_a=\inf\{t>0: X_t\le a\},
$$
then we have
\begin{equation}\label{ksv-e:firstexittime}
{\mathbb P}_x(\tau_a=\sigma_a)=1, \quad x>a.
\end{equation}

Let $G_{(0, \infty)}(x, y)$ be the Green function of
$X$ in $(0, \infty)$.
Then we have the following result.

\begin{prop}\label{ksv-p:Greenf4kpXonhalfline} For any $x, y>0$ we have
$$
G_{(0, \infty)}(x, y)=\left\{\begin{array}{ll}
\int^x_0v(z)v(y+z-x)dz, & x\le y,\\
\int^x_{x-y}v(z)v(y+z-x)dz, & x>y.
\end{array}\right.
$$
\end{prop}

\pf
Let $X^{(0,\infty)}$ be the process obtained by killing $X$ upon exiting from
$(0, \infty)$.
By using (\ref{ksv-e:firstexittime}) above and \cite[Theorem 20,
p.~176]{Ber} we get that for any nonnegative function on $f$ on $(0,
\infty)$,
\begin{equation}\label{ksv-e:e1inpfoformforgf}
{\mathbb E}_x\left[ \int_0^{\infty} f(X^{(0, \infty)}_t)\, dt\right]=k
\int^{\infty}_0 \int^x_0v(z)f(x+z-y)v(y) dz dy\, ,
\end{equation}
where $k$ is the constant depending on the normalization of the
local time of the process $X$ reflected at its supremum. We choose
$k=1$. Then
\begin{eqnarray}\label{ksv-e:e2inpfoformforgf}
&&{\mathbb E}_x\left[ \int_0^{\infty} f(X^{(0, \infty)}_t)\, dt\right]
\,=\,\int_0^{\infty} \, v(y)\int_0^x \, v(z)f(x+y-z) dz dy\nonumber\\
&&=\int_0^x \, v(z)\int_0^{\infty}  v(y)f(x+y-z) dy dz
\,=\,\int_0^x \, v(z)\int_{x-z}^{\infty}\, v(w+z-x)f(w) dw dz\nonumber\\
&&=\int_0^x f(w) \int_{x-w}^x \,  v(z)v(w+z-x) dzdw+ \int_x^{\infty}
f(w) \int_0^x \, v(z)v(w+z-x)  dzdw\, .\nonumber \\
&&
\end{eqnarray}
On the other hand,
\begin{equation}\label{ksv-e:e3inpfoformforgf}
{\mathbb E}_x\left[\int_0^{\infty} f(X^{(0, \infty)}_t)\, dt\right]
=\int_0^{\infty}G_{(0, \infty)}(x,w)f(w)\, dw.
\end{equation}
By comparing (\ref{ksv-e:e2inpfoformforgf}) and
(\ref{ksv-e:e3inpfoformforgf}) we arrive at our desired conclusion. \qed

For any $r>0$, let $G_{(0, r)}$ be the Green function
of $X$ in  $(0,
r)$.
Then we have the following result.

\begin{prop}\label{ksv-p:upbdongfofkpinfiniteinterval}
For all $r>0$ and all $x\in (0,r)$
$$
\int_0^r G_{(0,r)}(x,y)\, dy \le 2 V(x) V(r)\, .
$$
In particular, for any $R>0$, there exists $C_6=C_6(R)>0$ such that
for all $r\in (0,R)$ and all $x\in (0,r)$,
$$
\int^r_0 G_{(0, r)}(x, y)dy \le  C_6 (\phi(r^{-2})\phi(x^{-2}))^{-1/2}
\asymp\frac{r^{\alpha/2}}{(\ell(r^{-2}))^{1/2}} \frac{x^{\alpha/2}}{(\ell(x^{-2}))^{-1/2}}\, .
$$

\end{prop}

\pf For any $x\in (0, r)$, we have
\begin{eqnarray*}
&&\int^r_0G_{(0, r)}(x, y)dy
\le \int^r_0G_{(0, \infty)}(x, y)dy\\
&&=\int^x_0\int^x_{x-y}v(z)v(y+z-x)dzdy+
\int^r_x\int^x_0v(z)v(y+z-x)dzdy\\
&&=\int^x_0v(z)\int^x_{x-z}v(y+z-x)dydz
+\int^x_0v(z)\int^r_xv(y+z-x)dydz\, \le\, 2\,V(r)\,V(x).
\end{eqnarray*}
Now the desired conclusion follows easily from Proposition
\ref{ksv-p:abofgf4lhpat0}.
\qed

As a consequence of the result above, we immediately get the
following.

\begin{corollary}\label{ksv-p:upbdongfofkpinfiniteinterval2}
For all $r>0$ and all $x\in (0,r)$
$$
\int_0^r G_{(0,r)}(x,y)\, dy \le 2 V(r)\big(V(x) \wedge V(r-x)\big)\, .
$$
In particular, for any $R>0$, there exists $C_7=C_7(R)>0$ such that for all $x\in (0,
r)$, and  $r\in (0, R)$,
\begin{eqnarray*}
\int^r_0G_{(0, r)}(x, y)dy&\le&
C_7 (\phi(r^{-2}))^{-1/2}\left((\phi(x^{-2}))^{-1/2}\wedge
(\phi((r-x)^{-2}))^{-1/2}\right) \\
&\asymp&\frac{r^{\alpha/2}}{(\ell(r^{-2}))^{1/2}}
 \left(\frac{x^{\alpha/2}}{(\ell(x^{-2}))^{1/2}}\wedge
\frac{(r-x)^{\alpha/2}}{(\ell((r-x)^{-2}))^{1/2}}\right)\, .
\end{eqnarray*}
\end{corollary}
\pf The first inequality is a consequence of the identity $\int^r_0
G_{(0, r)}(x, y)dy=\int^r_0 G_{(0, r)}(r-x, y)dy$ which is true by
symmetry of the process $X$. The second one now follows exactly as
in the proof of Proposition \ref{ksv-p:upbdongfofkpinfiniteinterval}.
\qed

\begin{remark}\label{ksv-r:upbdongfofkpinfiniteinterval2}{\rm
With self-explanatory notation, an immediate consequence of the
above corollary is the following estimate
\begin{equation}\label{ksv-e:upbdongfofkpinfiniteinterval2}
\int_{-r}^r G_{(-r,r)}(x,y)\, dy \le 2 V(2r)\big(V(r+x) \wedge V(r-x)\big)\, .
\end{equation}
}
\end{remark}

\section{Harnack inequality and Boundary Harnack principle}\label{ksv-sec-hibhp}

From now on we will always assume that $X$ is a subordinate Brownian
motion in ${\mathbb R}^d$. Recall that {\bf (H)} is the standing assumptions
on the Laplace exponent $\phi$.
The goal of this section is to show that
the Harnack inequality and the boundary Harnack principle hold for
$X$. The infinitesimal generator ${\bf L}$ of the corresponding semigroup
is given by
\begin{equation}\label{ksv-3.1}
{\bf L} f(x)=\int_{{\mathbb R}^d}\left( f(x+y)-f(x)-y\cdot \nabla f(x)
{\bf 1}_{\{|y|\le1\}}
 \right)\, J(y)dy
\end{equation}
for $f\in C_b^2({\mathbb R}^d)$. Moreover, for every $f\in C_b^2({\mathbb R}^d)$
$$
f(X_t)-f(X_0)-\int_0^t {\bf L} f(X_s)\, ds
$$
is a ${\mathbb P}_x$-martingale for every $x\in {\mathbb R}^d$. We recall the L\'evy
system formula for $X$ which describes the jumps of the process $X$:
for any non-negative measurable function $f$ on ${\mathbb R}_+ \times
{\mathbb R}^d\times {\mathbb R}^d$ with $f(s, y, y)=0$ for all $y\in {\mathbb R}^d$, any
stopping time $T$ (with respect to the filtration of
$X$)
and any
$x\in {\mathbb R}^d$,
\begin{equation}\label{ksv-e:levy}
{\mathbb E}_x \left[\sum_{s\le T} f(s,X_{s-}, X_s) \right]= {\mathbb E}_x \left[
\int_0^T \left( \int_{{\mathbb R}^d} f(s,X_s, y) J(X_s,y) dy \right) ds
\right].
\end{equation}
 (See, for example, \cite[Proof of Lemma 4.7]{CK1} and \cite[Appendix A]{CK2}.)

\subsection{Harnack inequality}\label{ksv-ss:hi}

It follows from Theorem \ref{ksv-t:Jorigin} and the 0-version of
\cite[Propositions 1.5.8 and 1.5.10]{BGT} that
\begin{equation}\label{ksv-e:svphi1}
r^{-2}\int_0^r s^{d+1} j(s)
ds\asymp\frac{\ell(r^{-2})}{r^\alpha}\asymp \phi(r^{-2}), \qquad
r\to 0
\end{equation}
and
\begin{equation}\label{ksv-e:svphi2}
\int_r^{\infty} s^{d-1} j(s) ds
\asymp\frac{\ell(r^{-2})}{r^\alpha}\asymp \phi(r^{-2}), \qquad r\to
0.
\end{equation}

For any open set $D$, we use $\tau_D$ to denote the first exit time
from $D$, i.e., $\tau_D=\inf\{t>0: \, X_t\notin D\}$.

\begin{lemma}\label{ksv-L3.1} There exists a constant $C_8>0$ such that for
every $r\in (0,1)$ and every $t>0$,
$$
{\mathbb P}_x\left(\sup_{s\le t} |X_s-X_0|>r\right) \le C_8 \phi(r^{-2}) t\, .
$$
\end{lemma}
\pf It suffices to prove the lemma for $x=0$. Let $f\in
C^2_b({\mathbb R}^d)$, $0\leq f \leq 1$, $f(0)=0$, and $f(y)=1$ for all
$|y|\ge 1$. Let $c_1=\sup_{y}\sum_{j,k} |(\partial^2/\partial
y_j\partial y_k) f(y)|$. Then $|f(z+y)-f(z) -y\cdot \nabla f(z)|\le
\frac{c_1}{2} |y|^2$. For $r\in (0,1)$, let $f_r(y)=f(y/r)$. Then the
following estimate is valid:
\begin{eqnarray*}
|f_r(z+y)-f_r(z) -y\cdot \nabla f_r(z){\bf 1}_{\{|y|\le r\}}| &\le&
\frac{c_1}{2}
\frac{|y|^2}{r^2}{\bf 1}_{\{|y|\le r\}} +  {\bf 1}_{\{|y|\ge r\}}\\
 &\le& c_2({\bf 1}_{\{|y|\le r\}}
\frac{|y|^2}{r^2}+{\bf 1}_{\{|y|\ge r\}})\, .
\end{eqnarray*}
By using \eqref{ksv-e:svphi1} and \eqref{ksv-e:svphi2}, we get the following
estimate:
\begin{eqnarray}\label{ksv-referee0}
|{\bf L} f_r(z)| &\le & \int_{{\mathbb R}^d}
|f_r(z+y)-f_r(z) -y\cdot \nabla f_r(z){\bf 1}_{(|y|\le r)}| \, J(y)dy \nonumber \\
& \le & c_2\int_{{\mathbb R}^d} \left({\bf 1}_{\{|y|\le r\}}
\frac{|y|^2}{r^2}+{\bf 1}_{\{|y|\ge r\}}\right)\, J(y)dy\nonumber \\
& \le &
C_8\phi(r^{-2}) \, ,
\end{eqnarray}
where the constant $C_8$ is independent of $r$. Further, by the
martingale property,
\begin{equation}\label{ksv-referee}
{\mathbb E}_0 f_r(
X_{\tau_{B(0,r)}\wedge t}
) - f_r(0)= {\mathbb E}_0
\int_0^{\tau_{B(0,r)}\wedge t} {\bf L} f_r(X_s)\, ds
\end{equation}
implying the estimate
$$
{\mathbb E}_0 f_{r}(
X_{\tau_{B(0,r)}\wedge t}) \leq
C_8\phi(r^{-2}) t\, .
$$
If $X$ exits $B(0,r)$ before time $t$, then
$f_{r}(
X_{\tau_{B(0,r)}\wedge t})=1$, so the left hand side is larger
than ${\mathbb P}_0(\tau_{B(0,r)} \le t)$. \qed

\begin{lemma}\label{ksv-L3.2}
For every $r\in (0,1)$, and every $x\in {\mathbb R}^d$,
$$
\inf_{z\in B(x,r/2)} {\mathbb E}_z \left[\tau_{B(x,r)} \right] \geq
\frac{1}{C_8 \phi((r/2)^{-2})}\, ,
$$
where $C_8$ is the constant from Lemma \ref{ksv-L3.1}.
\end{lemma}
\pf
Using \eqref{ksv-referee0} and \eqref{ksv-referee}
we get that for any $t>0$ and $z\in B(x,r/2)$,
\begin{eqnarray*}
{\mathbb P}_0(\tau_{B(0,r/2)} \le t)&\le&
C_8\phi((r/2)^{-2}){\mathbb E}_0
\left[\tau_{B(0,r/2)}\wedge t\right]\\
&=&
C_8\phi((r/2)^{-2}){\mathbb E}_{z}
\left[\tau_{B(z,r/2)}\wedge t\right]\\
&\le &
C_8\phi((r/2)^{-2}){\mathbb E}_{z}
\left[\tau_{B(x,r)}\wedge t\right].
\end{eqnarray*}
Letting $t\to\infty$, we immediately
get the desired conclusion.
\qed

\begin{lemma}\label{ksv-L3.3}
There exists a constant $C_9>0$ such that for every $r\in (0,1)$ and
every $x\in {\mathbb R}^d$,
$$
\sup_{z\in B(x,r)}
{\mathbb E}_z \left[{\tau}_{B(x,r)}\right] \leq
\frac{C_9}{\phi(r^{-2})}\, .
$$
\end{lemma}
\pf
Let $r\in (0,1)$, and let $x\in {\mathbb R}^d$. Using the L\'evy system
formula \eqref{ksv-e:levy},  we get
\begin{eqnarray*}
1 & \geq & {\mathbb P}_z (|
X_{{\tau}_{B(x,r)}}-x|>r ) \\
 & = & \int_{B(x,r)}G_{B(x,r)}(z,y) \int_{
 \overline{B(x,r)}^c
 } j(|u-y|)\, du \, dy \, ,
\end{eqnarray*}
where $G_{B(x,r)}$ denotes the Green function of the process
$X$ in  $B(x,r)$.
Now we estimate the inner integral.
Let $y\in B(x,r)$, $u\in
 \overline{B(x,r)}^c$. If $u\in B(x,2)$, then $|u-y|\le
2|u-x|$, while for $u\notin B(x,2)$ we use $|u-y|\le |u-x|+1$. Then
\begin{eqnarray*}
\lefteqn{\int_{
 \overline{B(x,r)}^c} j(|u-y|)\, du}\\
  & = & \int_{
 \overline{B(x,r)}^c\cap B(x,2)}
j(|u-y|)\, du+
\int_{
 \overline{B(x,r)}^c\cap B(x,2)^c} j(|u-y|)\, du \\
& \ge &\int_{
 \overline{B(x,r)}^c\cap B(x,2)} j(2|u-x|)\, du+
\int_{
 \overline{B(x,r)}^c\cap B(x,2)^c} j(|u-x|+1)\, du \\
& \ge & \int_{
 \overline{B(x,r)}^c\cap B(x,2)}c^{-1} j(|u-x|)\, du+
\int_{
 \overline{B(x,r)}^c\cap B(x,2)^c} c^{-1}j(|u-x|)\, du \\
&=& \int_{
 \overline{B(x,r)}^c}c^{-1}j(|u-x|)\, du\, ,
\end{eqnarray*}
where in the next to last line we used  \eqref{ksv-H:1} and \eqref{ksv-H:2}.
Now, It follows from \eqref{ksv-e:svphi2} that
\begin{eqnarray*}
1 &\ge &  \int_{B(x,r)}G_{B(x,r)}(z,y)\, dy
\int_{
 \overline{B(x,r)}^c}c^{-1}j(|u-x|)\, du \\
& = &
{\mathbb E}_z \left[\tau_{B(x,r)} \right]
 c^{-1}\, c_1\int_r^{\infty} v^{d-1} j(v)\, dv \\
& = & c_2\phi(r^{-2})
{\mathbb E}_z \left[\tau_{B(x,r)} \right]
\end{eqnarray*}
which implies the lemma. \qed

An improved version of the above lemma will be given in Proposition \ref{ksv-l:tau} later
on.

\begin{lemma}\label{ksv-L3.4}
There exists a constant $C_{10}>0$ such that for every $r\in (0,1)$,
every $x\in {\mathbb R}^d$, and any $A\subset B(x,r)$
$$
{\mathbb P}_y\left(T_A < {\tau}_{B(x,3r)}\right) \geq C_{10} \frac{|A|}{|B(x, r)|}, \qquad
 \textrm{for all }y\in B(x,2r)\, .
$$
\end{lemma}
\pf Without loss of generality assume that ${\mathbb P}_y(T_A < {\tau}_{B(x,3r)})<1/4$.
Set $\tau={\tau}_{B(x,3r)}$. By Lemma \ref{ksv-L3.1},
${\mathbb P}_y(\tau\leq t) \leq {\mathbb P}_y(\tau_{B(y,r)}\leq t) \leq c_1
\phi(r^{-2}) t$. Choose $t_0= 1/(4c_1 \phi(r^{-2}))$, so that
${\mathbb P}_y(\tau\leq t_0) \leq 1/4$. Further, if $z\in B(x, 3r)$ and $u\in
A \subset B(x,r)$, then $|u-z| \leq 4r$.
 Since $j$ is decreasing,
$j(|u-z|) \geq j(4r)$. Thus,
\begin{eqnarray*}
{\mathbb P}_y (T_A < \tau) & \geq & {\mathbb E}_y \sum_{s\leq T_A \wedge \tau \wedge
t_0}
{\bf 1}_{\{X_{s-}\neq X_s, X_s\in A\}} \\
& = & {\mathbb E}_y \int_0^{T_A \wedge \tau \wedge t_0}
\int_A  j(|u-X_s|)\, du \, ds \\
& \geq & {\mathbb E}_y \int_0^{T_A \wedge \tau \wedge t_0} \int_A  j(4r)\,
du \, ds \\
& = &  j(4r) |A| {\mathbb E}_y[T_A \wedge \tau \wedge t_0] \, ,
\end{eqnarray*}
where in the second line we used properties of the L\'evy system.
Next,
\begin{eqnarray*}
{\mathbb E}_y[T_A\wedge \tau \wedge t_0] & \ge & {\mathbb E}_y[t_0; \, T_A\geq \tau \geq t_0] \\
& = & t_0 {\mathbb P}_y(T_A \ge \tau \ge t_0) \\
& \ge & t_0[1-{\mathbb P}_y(T_A < \tau)-{\mathbb P}_y(\tau <t_0)] \\
& \ge & \frac{t_0}{2} = \frac{1}{8 c_1 \phi(r^{-2})}\, .
\end{eqnarray*}
The last two displays give that
$$
{\mathbb P}_y (T_A < \tau) \geq  j(4r) |A| \frac{1}{8 c_1 \phi(r^{-2})} =
\frac{1}{8 c_1} |A| \frac{j(4r)}{\phi(r^{-2})}.
$$
The claim now follows immediately from \eqref{ksv-e:reg-var} and
Theorem \ref{ksv-t:Jorigin}. \qed

\begin{lemma}\label{ksv-L3.5}
There exist positive constant $C_{11}$ and $C_{12}$, such that if
$r\in (0,1)$, $x\in {\mathbb R}^d$, $z\in B(x,r)$, and $H$ is a bounded
nonnegative function with support in $B(x,2r)^c$, then
$$
{\mathbb E}_z H(
X_{{\tau}_{B(x,r)}}) \leq C_{11} {\mathbb E}_z [{\tau}_{B(x,r)}] \int H(y)
j(|y-x|) \, dy \, ,
$$
and
$$
{\mathbb E}_z H(
X_{{\tau}_{B(x,r)}}) \geq C_{12} {\mathbb E}_z [{\tau}_{B(x,r)}] \int H(y)
j(|y-x|) \, dy \, .
$$
\end{lemma}
\pf Let $y\in B(x,r)$ and $u\in B(x,2r)^c$. If $u\in B(x,2)$ we use
the estimates
\begin{equation}\label{ksv-3.7}
2^{-1}|u-x|\le |u-y| \le 2|u-x|,
\end{equation}
while if $u\notin B(x,2)$ we use
\begin{equation}\label{ksv-3.8}
|u-x|-1\le |u-y|\le |u-x|+1.
\end{equation}
Let $B\subset B(x,2r)^c$. Then using the L\'evy system we get
$$
{\mathbb E}_z \left[ {\bf 1}_B(X_{\tau_{B(x,r)}}) \right]  =  {\mathbb E}_z \int_0^{\tau_{B(x,r)}} \int_B
j(|u-X_s|)\, du\, ds\, .
$$
By use of \eqref{ksv-H:1}, \eqref{ksv-H:2}, \eqref{ksv-3.7}, and \eqref{ksv-3.8},
 the inner integral is estimated as follows:
\begin{eqnarray*}
\int_B j(|u-X_s|)\, du &=& \int_{B\cap B(x,2)} j(|u-X_s|)\, du +
\int_{B\cap B(x,2)^c} j(|u-X_s|)\, du  \\
& \le & \int_{B\cap B(x,2)} j(2^{-1}|u-x|)\, du +
\int_{B\cap B(x,2)^c} j(|u-x|-1)\, du \\
& \le &\int_{B\cap B(x,2)} c j(|u-x|)\, du +
\int_{B\cap B(x,2)^c} c j(|u-x|)\, du \\
& = & c \int_B j(|u-x|)\, du.
\end{eqnarray*}
Therefore
\begin{eqnarray*}
{\mathbb E}_z \left[ {\bf 1}_B(X_{\tau_{B(x,r)}}) \right] & \le & {\mathbb E}_z \int_0^{\tau_{B(x,r)}}
c \int_B j(|u-x|)\, du \\
& = & c\,  {\mathbb E}_z (\tau_{B(x,r)}) \int {\bf 1}_B(u) j(|u-x|)\, du\, .
\end{eqnarray*}
Using linearity we get the above inequality when ${\bf 1}_B$ is replaced
by a simple function. Approximating $H$ by simple functions and
taking limits we have the first inequality in the statement of the
lemma.

The second inequality is proved in the same way. \qed

\begin{defn}\label{ksv-def:har1}
Let $D$ be an open subset of ${\mathbb R}^d$. A function $u$ defined on
${\mathbb R}^d$ is said to be

\begin{description}
\item{(1)}  harmonic in $D$ with respect to $X$ if
$$
{\mathbb E}_x\left[|u(X_{\tau_{B}})|\right] <\infty \quad \hbox{ and } \quad
u(x)= {\mathbb E}_x\left[u(X_{\tau_{B}})\right], \qquad x\in B,
$$
for every open set $B$ whose closure is a compact subset of $D$;

\item{(2})
regular harmonic in $D$ with respect to $X$
if
it is harmonic in $D$ with respect to $X$ and
for each $x \in D$,
$u(x)= {\mathbb E}_x\left[u(X_{\tau_{D}})\right].$
\end{description}
\end{defn}

Now we give the proof of Harnack inequality.
The proof below is basically the proof given in \cite{SV04}
which is an adaptation of the proof given in \cite{BL02a}.
However, the proof below corrects some typos in the proof given in \cite{SV04}.

\begin{thm}\label{ksv-T:Har}
There exists $C_{13}>0$ such that, for any $r\in (0, 1/4)$,
$x_0\in{\mathbb R}^d$, and any function $u$ which is nonnegative on ${\mathbb R}^d$
and harmonic with respect to $X$ in $B(x_0,
17r)$, we have
$$
u(x)\le C_{13} u(y), \quad \textrm{for all }x, y\in B(x_0, r).
$$
\end{thm}
\pf
Without loss of generality we may assume that $u$ is strictly
positive in $B(x_0, 16r)$. Indeed, if $u(x)=0$ for some $x\in B(x_0,
16r)$, then by harmonicity $ 0=u(x)={\mathbb E}_x [u(X_{\tau_B})] $
for $x
\in B=B(x,\epsilon) \subset B(x_0, 16r)$. This and the fact that the
Levy measure of $X$ is supported on all of ${\mathbb R}^d$ and has a density
imply that $u=0$ a.e. with respect to Lebesgue measure. Moreover,
by the harmonicity,  for every $y \in B(x_0, 16r)$,
$u(y)={\mathbb E}_y[u(X_{\tau_B} )]=0$ where $B=B(y,\delta)\subset B(x_0,
16r)$. Therefore, if $u(x)=0$ for some $x$, then $u$ is identically
zero in $B(x_0, 16 r)$ and there is nothing to prove.

We first assume $u$ is bounded on ${\mathbb R}^d$.
Using the harmonicity of $u$ and Lemma
\ref{ksv-L3.4}, one can show that $u$ is bounded from below on $B(x_0,
r)$ by a positive number. To see this, let $\epsilon>0$ be such that
$F=\{x\in B(x_0, 3r)\setminus B(x_0, 2r): u(x)>\epsilon\}$ has
positive Lebesgue measure. Take a compact subset $K$ of $F$ so that
it has positive Lebesgue measure. Then by Lemma \ref{ksv-L3.4}, for
$x\in B(x_0, r)$, we have
$$
u(x)\,=\,{\mathbb E}_x\left[u(
X_{T_K\wedge \tau_{B(x_0,9r)}}
)
\right]
\,>\,c\,\epsilon\,\frac{|K|}{|B(x_0, 3r)|},
$$
for some $c>0$.  By taking a constant multiple of $u$ we may assume
that $\inf_{B(x_0, r)}u =1/2$. Choose $z_0\in B(x_0, r)$ such that
$u(z_0)\le 1$. We want to show that $u$ is bounded above in $B(x_0,
r)$ by a positive constant independent of $u$ and $r\in (0, 1/4)$.
We will establish this by contradiction: If there exists a point
$x\in B(x_0, r)$ with $
u(x)=K$ where $K$ is too large, we can obtain
a sequence of points in $B(x_0, 2r)$ along which $u$ is unbounded.

Using
Lemmas
\ref{ksv-L3.2}, \ref{ksv-L3.3} and  \ref{ksv-L3.5}, one can see that
there exists $c_1>0$ such that if $x\in {\mathbb R}^d$, $s\in (0, 1)$ and $H$
is nonnegative bounded function with support in $B(x, 2s)^c$, then
for any $y, z\in B(x, s/2)$,
\begin{equation}\label{ksv-e:2.1}
{\mathbb E}_z H(
X_{\tau_{B(x, s)}}
)\,\le \,c_1\,{\mathbb E}_y H(
X_{\tau_{B(x, s)}}
).
\end{equation}
By Lemma \ref{ksv-L3.4}, there exists $c_2>0$ such that if $A\subset
B(x_0, 4r)$ then
\begin{equation}\label{ksv-e:2.2}
{\mathbb P}_y\left(T_A<\tau_{B(x_0, 16r)}\right)\,\ge\, c_2\,\frac{|A|}{|B(x_0, 4r)|}, \quad
\forall y\in B(x_0, 8r).
\end{equation}
Again by Lemma \ref{ksv-L3.4}, there exists $c_3>0$ such that if
$x\in{\mathbb R}^d$, $s\in (0, 1)$ and $F\subset B(x, s/3)$ with $|F|/|B(x,
s/3)|\ge 1/3$, then
\begin{equation}\label{ksv-e:2.3}
{\mathbb P}_x\left(T_F<\tau_{B(x, s)}\right)\,\ge\, c_3.
\end{equation}
Let
\begin{equation}\label{ksv-e:2.4}
\eta=\frac{c_3}3,\,\,\,\,\,\,\,\,\,\,\,
\zeta=(\frac13\wedge\frac1{c_1})\eta.
\end{equation}
Now suppose there exists $x\in B(x_0, r)$ with $u(x)=K$ for
$K>K_0:=\frac{2|B(x_0, 1)|}{c_2\zeta}\vee\frac{2(12)^d}{c_2\zeta}$. Let
$s$ be chosen so that
\begin{equation}\label{ksv-e:2.5}
|B(x, \frac{s}3)|=\frac{2|B(x_0, 4r)|}{c_2\zeta K}<1.
\end{equation}
Note that this implies
\begin{equation}\label{ksv-e:2.6}
s=12\left(\frac2{c_2\zeta}\right)^{1/d}rK^{-1/d}<r.
\end{equation}
Let us write $B_s$ for $B(x, s)$, $\tau_s$ for $\tau_{B(x, s)}$, and
similarly for $B_{2s}$ and $\tau_{2s}$. Let $A$ be a compact subset
of
$$
A'=\{y\in B(x, \frac{s}3): u(y)\ge \zeta K\}.
$$
It is well known that $u(X_t)$ is right continuous in
$[0,\tau_{B(x_0, 16r)})$. Since $z_0\in B(x_0, r)$ and $A'\subset B(x,
\frac{s}3)\subset B(x_0, 2r)$, we can apply (\ref{ksv-e:2.2}) to get
\begin{eqnarray*}
1&\ge&u(z_0)\ge {\mathbb E}_{z_0}[u(
X_{T_A\wedge\tau_{B(x_0, 16r)}}
){\bf 1}_{\{T_A<
\tau_{B(x_0, 16r)}\}}]\\
&\ge&\zeta K{\mathbb P}_{z_0}(T_A<\tau_{B(x_0, 16r)})\\
&\ge&c_2\zeta K\frac{|A|}{|B(x_0, 4r)|}.
\end{eqnarray*}
Hence
$$
\frac{|A|}{|B(x, \frac{s}3)|}\le\frac{|B(x_0, 4r)|} {c_2\zeta K
|B(x, \frac{s}3)|}=\frac12.
$$
This implies that $|A'|/|B(x, s/3)|\le 1/2$. Let $F$ be a compact
subset of $B(x, s/3)\setminus A'$ such that
\begin{equation}\label{ksv-e:2.7}
\frac{|F|}{|B(x, \frac{s}3)|}\ge \frac13.
\end{equation}
Let $H=u\cdot{\bf 1}_{B_{2s}^c}$. We claim that
$$
{\mathbb E}_x[u(
X_{\tau_s}
);
X_{\tau_s}
\notin B_{2s}]\le\eta K.
$$
If not, ${\mathbb E}_x H(
X_{\tau_s})>\eta K$, and by (\ref{ksv-e:2.1}), for all
$y\in B(x, s/3)$, we have
\begin{eqnarray*}
u(y)&=&{\mathbb E}_y u(
X_{\tau_s})\ge {\mathbb E}_y[u(
X_{\tau_s});
X_{\tau_s}\notin B_{2s}]\\
&\ge& c_1^{-1}{\mathbb E}_x H(
X_{\tau_s})\ge c_1^{-1}\eta K\ge \zeta K,
\end{eqnarray*}
contradicting (\ref{ksv-e:2.7}) and the definition of $A'$.

Let $M=\sup_{B_{2s}}u$. We then have
\begin{eqnarray*}
K&=&u(x)=
{\mathbb E}_x [u(
X_{\tau_s  \wedge T_F}
)]\\
&=&{\mathbb E}_x[u(
X_{T_F}
); T_F<\tau_s]+
{\mathbb E}_x[u(
X_{\tau_s}); \tau_s<T_F,
X_{\tau_s}\in B_{2s}]\\
&&\,\,\,+
{\mathbb E}_x[u(
X_{\tau_s}); \tau_s<T_F,
X_{\tau_s}\notin B_{2s}]\\
&\le& \zeta K{\mathbb P}_x(T_F<\tau_s)+M{\mathbb P}_x(\tau_s<T_F)+\eta K\\
&=&\zeta K{\mathbb P}_x(T_F<\tau_s)+M(1-{\mathbb P}_x(T_F<\tau_s))+\eta K,
\end{eqnarray*}
or equivalently
$$
\frac{M}{K}\ge\frac{1-\eta-\zeta}{1-{\mathbb P}_x(T_F<\tau_s)} +\zeta .
$$
Using (\ref{ksv-e:2.3}) and (\ref{ksv-e:2.4}) we see that there exists
$\beta>0$ such that $M\ge K(1+2\beta)$. Therefore there exists
$x'\in B(x, 2s)$ with $u(x')\ge K(1+\beta)$.

Now suppose there exists $x_1\in B(x_0, r)$ with $u(x_1)=K_1>K_0$.
Define $s_1$ in terms of $K_1$ analogously to (\ref{ksv-e:2.5}). Using
the above argument (with $x_1$ replacing $x$ and $x_2$ replacing
$x'$), there exists $x_2\in B(x_1, 2s_1)$ with $u(x_2)=K_2\ge
(1+\beta)K_1$. We continue and obtain $s_2$ and then $x_3$, $K_3$,
$s_3$, etc. Note that $x_{i+1}\in B(x_i, 2s_i)$ and $K_i\ge
(1+\beta)^{i-1}K_1$. In view of (\ref{ksv-e:2.6}),
\begin{eqnarray*}\sum_{i=0}^{\infty}
|x_{i+1}-x_i|&\le& r+ 24 \left(\frac2{c_2\zeta}\right)^{1/d}r
\sum_{i=1}^{\infty}K_i^{-1/d}\\
& \le & r + 24 \left(\frac2{c_2\zeta}\right)^{1/d}  K_1^{-1/d}r
\sum_{i=1}^{\infty} (1+\beta)^{-(i-1)/d}\\
&=& r + 24r \left(\frac2{c_2\zeta}\right)^{1/d}  K_1^{-1/d}r\sum^\infty_{i=0}(1+\beta)^{-i/d}\\
&=& r+ c_4rK^{-1/d}_1
\end{eqnarray*}
where $c_4:=24 (\frac2{c_2\zeta})^{1/d}\sum^\infty_{i=0}(1+\beta)^{-i/d}$.
So if $K_1>c^d_4 \vee K_0$ then we have a sequence $x_1, x_2, \dots$
contained in $B(x_0, 2r)$ with $u(x_i)\ge
(1+\beta)^{i-1}K_1\rightarrow\infty$, a contradiction to $u$ being
bounded. Therefore we can not take $K_1$ larger than $c^d_4\vee
K_0$, and thus $\sup_{y\in B(x_0, r)}u(y)\le c^d_4\vee K_0$, which
is what we set out to prove.

In the case that $u$ is unbounded, one can follow the simple limit
argument in the proof of \cite[Theorem 2.4]{SV04} to finish the
proof. \qed

By using the standard chain argument one can derive the following
form of Harnack inequality.

\begin{corollary}\label{ksv-c:hi}
For every $a \in (0,1)$, there exists $C_{14}=C_{14}(a)>0$ such that
for every $r \in (0, 1/4)$, $x_0 \in {\mathbb R}^d$, and any function $u$ which is nonnegative on ${\mathbb R}^d$
and harmonic with respect to $X$ in $B(x_0, r)$, we have
$$
u(x)\le C_{14} u(y), \quad \textrm{for all }x, y\in B(x_0, ar)\, .
$$
\end{corollary}

\subsection{Some estimates for the Poisson kernel}

Recall that for any open set $D$ in ${\mathbb R}^d$,
$\tau_D$
is the first exit time of $X$ from
$D$.

We recall  from Subsection \ref{ksv-ss:sbm} that $X$ has a transition
density $q(t, x, y)$, which is jointly continuous. Using this and the strong Markov property, one can easily check that
$$
q_D(t, x, y):=q(t, x, y)-{\mathbb E}_x[t>\tau_D, q(t-\tau_D, X_{\tau_D}, y)], \quad x, y \in D
$$
is continuous and the transition density of $X^D$.
 For any bounded open set $D\subset {\mathbb R}^d$, we
will use $G_D$ to denote the Green function of $X^D$, i.e.,
$$
G_D(x, y):=\int^\infty_0 q_D(t, x, y)dt, \quad x, y\in D.
$$
Note that $G_D(x,y)$ continuous  in $(D\times D)\setminus\{(x, x): x\in D\}$.
We will frequently use the well-known fact that
$G_D(\cdot, y)$ is harmonic in $D\setminus\{y\}$, and regular
harmonic in $D\setminus \overline{B(y,\varepsilon)}$ for every
$\varepsilon >0$.

Using the L\'{e}vy system for $X$, we know that for every bounded
open subset $D$, every $f \ge 0$ and all $x \in D$,
\begin{equation}\label{ksv-newls}
{\mathbb E}_x\left[f(X_{\tau_D});\,X_{\tau_D-} \not= X_{\tau_D}  \right] =
\int_{\overline{D}^c} \int_{D} G_D(x,z) J(z-y) dz f(y)dy.
\end{equation}
For notational convenience, we define
\begin{equation}\label{ksv-PK}
K_D(x,y)\,:=  \int_{D} G_D(x,z) J(z-y) dz, \qquad (x,y) \in D \times
\overline{D}^c.
\end{equation}
Thus \eqref{ksv-newls} can be simply written as
\begin{equation}\label{ksv-newls-2}
{\mathbb E}_x\left[f(X_{\tau_D});\,X_{\tau_D-} \not= X_{\tau_D}  \right]
=\int_{\overline{D}^c} K_D(x,y)f(y)dy\, ,
\end{equation}
revealing $K_D(x,y)$ as a density of the exit distribution of $X$
from $D$. The function $K_D(x,y)$ is called
the Poisson kernel of $X$. Using the continuity of $G_D$ and $J$, one can
easily check that $K_D$ is continuous on $D \times \overline{D}^c$.

The following proposition is an improvement of Lemma \ref{ksv-L3.3}. The
idea of the proof  comes from \cite{Sz2}.
\begin{prop}\label{ksv-l:tau}
For all $r>0$ and all $x_0\in {\mathbb R}^d$,
$$
{\mathbb E}_x[\tau_{B(x_0,r)}]\,\le\, 2V(2r) V(r-|x-x_0|)\, ,\qquad x\in B(x_0, r)\, .
$$
In particular, for any $R>0$, $r\in (0, R)$ and $x_0 \in {\mathbb R}^d$,
\begin{eqnarray*}
{\mathbb E}_x[\tau_{B(x_0,r)}]&\le & C_{7}\, (\phi(r^{-2})\phi((r-|x-x_0|)^{-2}))^{-1/2}\\
&\asymp &\frac{r^{\alpha/2}}{(\ell(r^{-2}))^{1/2}}\frac{(r-|x-x_0|)^{\alpha/2}}{(\ell((r-|x-x_0|)^{-2}))^{1/2}},
\qquad x\in B(x_0, r)\, ,
\end{eqnarray*}
where $C_7=C_7(R)$ is the constant form Proposition
\ref{ksv-p:upbdongfofkpinfiniteinterval2}.
\end{prop}

\pf Without loss of generality, we may assume that $x_0=0$. For
$x\neq 0$, put $Z_t=\frac{X_t\cdot x}{|x|}$. Then $Z_t$ is a L\'evy
process on ${\mathbb R}$ with
$$
{\mathbb E}(e^{i\theta Z_t})={\mathbb E}(e^{i\theta\frac{x}{|x|}\cdot X_t})
=e^{-t \phi(|\theta\frac{x}{|x|}|^2)}=e^{-t \phi(\theta^2)}
\qquad \theta\in {\mathbb R}.
$$
Thus $Z_t$ is of the type of one-dimensional subordinate Brownian
motion studied in Section \ref{ksv-ss:1dsbm}. It is easy to see that, if
$X_t\in B(0, r)$, then $|Z_t|<r$, hence
$$
{\mathbb E}_x[\tau_{B(0, r)}]\le {\mathbb E}_{|x|}[\tilde \tau],
$$
where $\tilde \tau=\inf\{t>0: |Z_t|\ge r\}$. Now the desired
conclusion follows easily from Proposition
\ref{ksv-p:upbdongfofkpinfiniteinterval2} (more precisely, from
\eqref{ksv-e:upbdongfofkpinfiniteinterval2}). \qed

As a consequence of Lemma \ref{ksv-L3.2}, Proposition \ref{ksv-l:tau} and
\eqref{ksv-PK}, we get the following result.

\begin{prop}\label{ksv-p:Poisson1}
There exist  $C_{15}, C_{16}>0$ such that for every $r \in (0, 1)$ and
$x_0 \in {\mathbb R}^d$,
\begin{eqnarray}\label{ksv-P1}
K_{B(x_0,r)}(x,y) &\le & C_{15} \, j(|y-x_0|-r) \big(\phi(r^{-2})\phi((r-|x-x_0|)^{-2})\big)^{-1/2}\\
&\asymp & j(|y-x_0|-r) \frac{r^{\alpha/2}}
{(\ell(r^{-2}))^{1/2}}\frac{(r-|x-x_0|)^{\alpha/2}}
{(\ell((r-|x-x_0|)^{-2}))^{1/2}}\, ,\nonumber
\end{eqnarray}
for all $(x,y) \in B(x_0,r)\times \overline{B(x_0,r)}^c$ and
\begin{equation}\label{ksv-P2}
K_{B(x_0,r)}(x_0,y) \,\ge\, C_{16}\,\frac{j(|y-x_0|)}{\phi((r/2)^{-2})}\asymp
j(|y-x_0|)\frac{r^\alpha}{\ell(r^{-2})}
\end{equation}
for all $y \in\overline{B(x_0,r)}^c$.
\end{prop}
\pf Without loss of generality, we assume $x_0=0$. For $z \in B(0,
r)$ and $r<|y|<2$
$$
|y|-r \le |y|-|z| \le |z-y| \le |z|+|y| \le r +|y| \le 2|y|
,
$$
and for $z \in B(0, r)$ and $y \in B(0, 2)^c$,
$$
|y|-r \le |y|-|z| \le |z-y| \le |z|+|y| \le r +|y|\le |y|+1.
$$
Thus by the monotonicity of $j$, \eqref{ksv-H:1} and \eqref{ksv-H:2}, there
exists a constant $c>0$ such that
$$
c j(|y|) \,\le\, j(|z-y|) \, \le \, j(|y|-r)\, ,    \qquad (z,y) \in
B(0,r) \times \overline{B(0,r)}^c.
$$
Applying the above inequality, Lemma \ref{ksv-L3.2} and Proposition
\ref{ksv-l:tau} to \eqref{ksv-PK}, we have proved the proposition. \qed

\begin{prop}\label{ksv-p:Poisson2}
For every $a \in (0,1)$,  $r \in (0, 1/4)$, $x_0 \in {\mathbb R}^d$ and $x_1, x_2 \in B(x_0, ar)$,
$$
K_{B(x_0,r)}(x_1,y) \,\le\, C_{14}  K_{B(x_0,r)}(x_2,y), \qquad y \in
\overline{B(x_0,r)}^{\, c}\, ,
$$
where $C_{14}=C_{14}(a)$ is the constant from Corollary \ref{ksv-c:hi}.
\end{prop}

\pf Let $a\in (0,1)$, $r\in (0,1/4)$ and $x_0\in {\mathbb R}^d$ be fixed.
For every Borel set $A\subset \overline{B(x_0,r)}^{\, c}$,
the function $x\mapsto {\mathbb P}_x(X_{\tau_{B(x_0,r)}}\in A)$ is
harmonic in $B(x_0,r)$. By Corollary \ref{ksv-c:hi} and \eqref{ksv-newls-2},
we have for all $x_1, x_2 \in B(x_0, ar)$,
\begin{eqnarray*}
\int_A K_{B(x_0,r)}(x_1,y)\, dy&=& {\mathbb P}_{x_1}(X_{\tau_{B(x_0,r)}}\in A)\\
&\le & C_{14} {\mathbb P}_{x_2}(X_{\tau_{B(x_0,r)}}\in A)=\int_A K_{B(x_0,r)}(x_2,y)\, dy\, .
\end{eqnarray*}
This implies that $K_{B(x_0,r)}(x_1,y)\le C_{14} K_{B(x_0,r)}(x_2,y)$
for a.e.~$y\in \overline{B(x_0,r)}^{\, c}$, and hence by
the continuity of $ K_{B(x_0,r)}(x,\cdot)$ for every
$y\in \overline{B(x_0,r)}^{\, c}$.\qed

The next inequalities will be used several times in the remainder of
this paper.

\begin{lemma}\label{ksv-l:l}
There exists $C>0$
such that
\begin{equation}\label{ksv-el1}
\frac{s^{\alpha/2}}{\left(\ell(s^{-2})\right)^{1/2}} \,\le \, C \,
\frac{r^{\alpha/2}}{\left(\ell(r^{-2})\right)^{1/2}}, \qquad
 0<s<r\le 4,
\end{equation}
\begin{equation}\label{ksv-el2}
\frac{s^{1-\alpha/2}}{\left(\ell(s^{-2})\right)^{1/2}} \,\le \, C \,
\frac{r^{1-\alpha/2}}{\left(\ell(r^{-2})\right)^{1/2}}, \qquad
 0<s<r\le 4,
\end{equation}
\begin{equation}\label{ksv-el7}
s^{1-\alpha/2} \,{\left(\ell(s^{-2})\right)^{1/2}} \,\le \, C \,
r^{1-\alpha/2}\,{\left(\ell(r^{-2})\right)^{1/2}}, \qquad
0<s<r\le 4,
\end{equation}
\begin{equation}\label{ksv-el3}
\int^{\infty}_r
\frac{\left(\ell(s^{-2})\right)^{1/2}}{s^{1+\alpha/2}}ds
 \,\le \, C \,
\frac{\left(\ell(r^{-2})\right)^{1/2}}{r^{\alpha/2}}, \qquad
 0<r\le 4,
\end{equation}
\begin{equation}\label{ksv-el6}
\int^{r}_0 \frac{\left(\ell(s^{-2})\right)^{1/2}}{s^{\alpha/2}}ds
 \,\le \, C \,
\frac{\left(\ell(r^{-2})\right)^{1/2}}{r^{\alpha/2-1}}, \qquad
 0<r\le 4,
\end{equation}

\begin{equation}\label{ksv-el4}
\int^{\infty}_r \frac{\ell(s^{-2})}{s^{1+\alpha}}ds
 \,\le \, C \,
\frac{\ell(r^{-2})}{r^{\alpha}}, \qquad   0<r\le 4,
\end{equation}
\begin{equation}\label{ksv-el5}
\int_{0}^r \frac{\ell(s^{-2})}{s^{\alpha-1}}ds
 \,\le \, C \,
\frac{\ell(r^{-2})}{r^{\alpha-2}}, \qquad   0<r\le 4,
\end{equation}
and
\begin{equation}\label{ksv-el8}
\int_{0}^r \frac{s^{\alpha-1}}{\ell(s^{-2})}ds
 \,\le \, C \,
\frac{r^{\alpha}}{\ell(r^{-2})}, \qquad \ 0<r\le 4.
\end{equation}
\end{lemma}
\pf The first three inequalities follow easily from \cite[Theorem
1.5.3]{BGT}, while the last five from the 0-version of
\cite[1.5.11]{BGT}. \qed

\begin{prop}\label{ksv-p:Poisson3}
For every $a \in (0,1)$, there exists $C_{17}=C_{17}(a)>0$ such that for every
$r \in (0, 1)$ and $x_0 \in {\mathbb R}^d$,
\begin{eqnarray*}
K_{B(x_0,r)}(x,y) \,&\le &\, C_{17}\,\frac{r^{\alpha/2-d}}{(\ell(r^{-2}))^{1/2}}
\frac{(\ell((|y-x_0|-r)^{-2}))^{1/2}} {(
|y-x_0|-r)^{\alpha/2}}\, ,\\
& & \qquad \qquad \qquad \forall x\in  B(x_0, ar),\,  y \in
\{r<|x_0-y| \le 2r\}\, .
\end{eqnarray*}
\end{prop}

\pf By Proposition \ref{ksv-p:Poisson2},
$$
K_{B(x_0,r)}(x,y) \le \frac{c_1}{r^d} \int_{B(x_0, a r)}
K_{B(x_0,r)}(w,y) dw
$$
for some constant $c_1=c_1(a)>0$. Thus from Proposition \ref{ksv-l:tau}, (\ref{ksv-P1}) and Remark \ref{ksv-r:abofgf4lhpat0} we have that
\begin{eqnarray*}
K_{B(x_0,r)}(x,y)&\le& \frac{
c_1}{r^d}\int_{B(x_0, r)}\int_{B(x_0,r)} G_{B(x_0,r)}(w,z)J(z-y) dz dw \\
&=& \frac{c_1}{r^d}\int_{B(x_0, r)} {\mathbb E}_z[\tau_{B(x_0,r)}]J(z-y) dz\\
& \le& \frac{c_2}{r^d}  \frac{r^{\alpha/2}}{(\ell(r^{-2}))^{1/2}}
\int_{B(x_0, r)}\frac{(r-|z-x_0|)^{\alpha/2}}{(\ell((r-|z-x_0|)^{-2}))^{1/2}} J(z-y)dz
\end{eqnarray*}
for some constant $c_2=c_2(a)>0$.
Now applying Theorem \ref{ksv-t:Jorigin}, we get
$$
K_{B(x_0,r)}(x,y) \, \le\,
\frac{c_3 r^{\alpha/2-d}}{(\ell(r^{-2}))^{1/2}} \int_{B(x_0,
r)}\frac{(r-|z-x_0|)^{\alpha/2}} {(\ell((r-|z-x_0|)^{-2}))^{1/2}}
\frac{\ell(|z-y|^{-2})} {|z-y|^{d+\alpha}}dz
$$
for some constant $c_3=c_3(a)>0$. Since $r-|z-x_0| \le |y-z| \le 3r
\le
3$,  from \eqref{ksv-el1} we see that
$$
\frac{(r-|z-x_0|)^{\alpha/2}} {(\ell((r-|z-x_0|)^{-2}))^{1/2}} \,
\le\, c_4
 \frac{(|y-z|)^{\alpha/2}}
{(\ell(|y-z|^{-2}))^{1/2}}
$$
for some constant $c_4>0$. Thus we have
\begin{eqnarray*}
K_{B(x_0,r)}(x,y) & \le&
\frac{c_5 r^{\alpha/2-d}}{(\ell(r^{-2}))^{1/2}} \int_{B(x_0, r)}
\frac{(\ell(|z-y|^{-2}))^{1/2}}{|z-y|^{d+\alpha/2}}dz\\
& \le&   \frac{c_5 r^{\alpha/2-d}}{(\ell(r^{-2}))^{1/2}} \int_{B(y,
|y-x_0|-r)^c}
\frac{(\ell(|z-y|^{-2}))^{1/2}}{|z-y|^{d+\alpha/2}}dz\\
&\le & \frac{c_6 r^{\alpha/2-d}}{(\ell(r^{-2}))^{1/2}}
\int_{|y-x_0|-r}^{\infty} \frac{\left(\ell(s^{-2})\right)^{1/2}}
{s^{1+\alpha/2}}ds
\end{eqnarray*}
for some constants $c_5=c_5(a)>0$ and  $c_6=c_6(a)>0$. Using
\eqref{ksv-el3} in the above equation, we conclude that
$$
K_{B(x_0,r)}(x,y) \,\le \,
\frac{c_7 r^{\alpha/2-d}}{(\ell(r^{-2}))^{1/2}}
\frac{(\ell((|y-x_0|-r)^{-2}))^{1/2}} {( |y-x_0|-r)^{\alpha/2}}
$$
for some constant $c_7=c_7(a)>0$. \qed

\begin{remark}\label{ksv-r:Poisson3}{\rm
Note that the right-hand side of the estimate can be replaced by $\frac{V(r)}{r^d V(|y-x_0|-r)}$.
}
\end{remark}

\subsection{Boundary Harnack principle}\label{ksv-ss:bhp}

In this subsection, we additionally assume that $\alpha\in (0, 2\wedge d)$ and
in the case $d\le 2$, we further assume \eqref{ksv-e:ass4trans}.

The proof of the boundary Harnack principle is basically
the proof given in \cite{KSV1}, which is adapted from \cite{Bog97, SW99}.
The following result is a generalization of \cite[Lemma 3.3]{SW99}.

\begin{lemma}\label{ksv-l2.1}
For every $a \in (0, 1)$, there exists a positive constant $C_{19}=C_{19}(a)>0$
 such that for any $r\in (0, 1)$ and any open
set $D$ with $D\subset B(0, r)$ we have
$$
{\mathbb P}_x\left(X_{\tau_D} \in B(0, r)^c\right) \,\le\, C_{19}\,r^{-\alpha}\,
\ell(r^{-2})\int_D G_D(x,y)dy, \qquad x \in D\cap B(0,
ar).
$$
\end{lemma}

\pf
We will use
$C^{\infty}_c({\mathbb R}^d)$ to denote
the space of infinitely differentiable functions with compact supports.
Recall that  ${\bf L}$ is the $L_2$-generator of $X$ in \eqref{ksv-3.1} and
that  $G(x,y)$ and $G_D(x,y)$ are the Green
functions of $X$
in ${\mathbb R}^d$ and $D$
 respectively. We have ${\bf L} \,
G(x,y)=-\delta_x(y)$ in the weak sense.  Since $ G_D(x,y)=G(x,y)
-{\mathbb E}_x[G(X_{\tau_D},y)] $,
we have, by the symmetry of ${\bf L}$, for any $x\in D$ and any
nonnegative $\phi \in C^{\infty}_c({\mathbb R}^d)$,
\begin{eqnarray*}
&&\int_D   G_D(x,y) {\bf L}  \phi(y)dy
=\int_{{\mathbb R}^d} G_D(x,y) {\bf L}  \phi(y)dy\\
 &&=
\int_{{\mathbb R}^d} G(x,y)  {\bf L}  \phi(y)dy-
 \int_{{\mathbb R}^d} {\mathbb E}_x[G(X_{\tau_D},y)]
{\bf L}  \phi(y)dy\\
 &&=
\int_{{\mathbb R}^d} G(x,y) {\bf L} \phi(y)dy- \int_{D^c}
\int_{{\mathbb R}^d} G(z,y) {\bf L} \phi(y)dy
{\mathbb P}_x(X_{\tau_D} \in dz)\\
 &&=-\phi(x)+ \int_{D^c}
\phi(z){\mathbb P}_x(X_{\tau_D} \in dz)
\,=\,-\phi(x)+{\mathbb E}_x[\phi(X_{\tau_D})].
\end{eqnarray*}
In particular, if $\phi(x)=0$ for $x\in D$, we have
\begin{equation}\label{ksv-har_gen}
{\mathbb E}_x\left[ \phi(X_{\tau_D})\right] = \int_D
G_D(x,y)  {\bf L}  \phi(y)dy.
\end{equation}

For fixed $a \in (0,1)$, take a
sequence of radial functions $\phi_m$ in $C^{\infty}_c({\mathbb R}^d)$ such
that $0\le \phi_m\le 1$,
\[
\phi_m(y)=\left\{
\begin{array}{lll}
0, & |y|< a\\
1, & 1\le |y|\le m+1\\
0, & |y|>m+2,
\end{array}
\right.
\]
and that $\sum_{i, j}|\frac{\partial^2}{\partial y_i\partial y_j}
\phi_m|$ is uniformly bounded.
Define $\phi_{m, r}(y)=\phi_m(\frac{y}{r})$ so that
$0\le \phi_{m, r}\le 1$,
 \begin{equation}\label{ksv-e:2.11}\phi_{m, r}(y)=
\begin{cases}
0, & |y|<ar\\
1, & r\le |y|\le r(m+1)\\
0, & |y|>r(m+2),
\end{cases}
 \quad \text{and} \quad  \sup_{y\in {\mathbb R}^d} \sum_{i,
j}\left|\frac{\partial^2}{\partial y_i\partial y_j} \phi_{m,
r}(y)\right| \,<\, c_1\, r^{-2}.
\end{equation}
We claim that there exists a constant $c_1=c_1(a)>0$ such that for
all $r\in (0, 1)$,
\begin{equation}\label{ksv-e2.1}
\sup_{m \ge 1}  \sup_{y\in {\mathbb R}^d} |{\bf L}\phi_{m,r}(y)|\,\le\,
c_1 r^{-\alpha} \, \ell(r^{-2}).
\end{equation}
In fact, by Proposition \ref{ksv-t:Jorigin} we have
\begin{eqnarray*}
&&  \left|\int_{{\mathbb R}^d} (\phi_{m,r}(x+y)-\phi_{m,r}(x)-(\nabla
\phi_{m,r}(x) \cdot y)1_{B(0, r)}(y))J(y)
 dy \right|\\
&&\le \left|\int_{\{|y|\le r\}}
(\phi_{m,r}(x+y)-\phi_{m,r}(x)-(\nabla
\phi_{m,r}(x) \cdot y)1_{B(0, r)}(y))J(y) dy\right|\\
&& \quad +2\int_{\{r<|y|\}}J(y) dy
\\
&&\le \frac{c_2}{r^2}\int_{\{|y|\le r \}}
|y|^2 J(y)dy
+2\int_{\{r<|y|\}}J(y) dy \\
&&\le \frac{c_3}{r^2}\int_{\{|y|\le r \}}
\frac1{|y|^{d+\alpha-2}} \ell(|y|^{-2})dy
\,+\,c_3\int_{\{r<|y|\}}   \frac1{|y|^{d+\alpha}} \ell(|y|^{-2})    dy
\\
&& \le
 \frac{c_4}{r^2}
 \int_{0}^r \frac{\ell(s^{-2})}{s^{\alpha-1}}ds\,+\,
c_4\int^{\infty}_r \frac{\ell(s^{-2})}{s^{1+\alpha}}ds.
\end{eqnarray*}
Applying \eqref{ksv-el4}-\eqref{ksv-el5} to the above equation, we get
$$
\left|\int_{{\mathbb R}^d} (\phi_{m,r}(x+y)-\phi_{m,r}(x)-(\nabla
\phi_{m,r}(x) \cdot y)1_{B(0, r)}(y))J(y)
 dy \right| \,\le \,
c_5\, r^{-\alpha}\, \ell(r^{-2}),$$
for some constant
$c_5=c_5(d, \alpha, \ell)>0$.
 So the claim follows.
Let  $A(x, a,b):=\{ y \in {\mathbb R}^d: a \le |y-x| <b  \}.$
When $D \subset B(0,r)$ for some $r\in (0, 1)$, we get,
by combining (\ref{ksv-har_gen})
and (\ref{ksv-e2.1}), that for any $x\in D\cap B(0, ar)$,
\begin{eqnarray*}
{\mathbb P}_x\left(X_{\tau_D} \in B(0, r)^c\right)\,&=&\, \lim_{m\to
\infty}{\mathbb P}_x\left(X_{\tau_D} \in A(0, r, (m+1)r)\right)\\
 \,&\le &\,
C\,r^{-\alpha} \, \ell(r^{-2})\int_D  G_D(x,y)dy.
\end{eqnarray*}
\qed

\begin{lemma}\label{ksv-l2.1_1}
There exists $C_{20}>0$ such that for any open set $D$ with  $B(A, \kappa r)\subset D\subset B(0, r)$
for some $r\in (0, 1)$ and
$\kappa\in (0, 1)$, we have that for every $x \in D \setminus B(A, \kappa r)$,
\begin{eqnarray*}
\lefteqn{\int_{D} G_D(x,y)   dy }\\
&  \le & C_{20}\, r^{\alpha}
\,\kappa^{-d-\alpha/2}\, \frac1{\ell((4r)^{-2})}\left(1+
\frac{\ell((\frac{\kappa r}{2})^{-2})}{\ell((4r)^{-2})}\right)
{\mathbb P}_x\left(X_{\tau_{D\setminus B(A, \kappa r)}} \in B(A, \kappa
r)\right).
\end{eqnarray*}
\end{lemma}

\pf Fix a point $x\in D\setminus B(A, \kappa r)$ and let $B:=B(A,
\frac{\kappa r}2)$. Note that, by the harmonicity of
$G_D(x,\,\cdot\,)$ in $D\setminus \{x\}$ with respect to $X$, we
have
\[
G_D(x,A)\,\ge\,\int_{D\cap \overline{B}^c}K_B(A, y)G_D(x,y)dy
\,\ge\,\int_{D\cap B(A, \frac{3\kappa r}4)^c}K_B(A, y)G_D(x,y)dy.
\]

Since $\frac{3\kappa r}4\le |y-A|\le 2r$ for $y\in B(A,
\frac{3\kappa r}4)^c\cap D$ and $j$ is a decreasing function, it
follows from \eqref{ksv-P2} in Proposition \ref{ksv-p:Poisson1} and Theorem \ref{ksv-t:Jorigin} that
\begin{eqnarray*}
G_D(x,A) &\ge& c_1\, \frac{(\frac{\kappa r}{2})^\alpha}{\ell\left((
\frac{\kappa r}{2})^{-2}\right)}\int_{D \cap B(A, \frac{3\kappa
r}4)^c}G_D(x,y)J(y-A) dy\\
&\ge& c_1\, j(2r)\, \frac{(\frac{\kappa r}{2})^\alpha}{\ell\left((
\frac{\kappa r}{2})^{-2}\right)}\int_{D \cap B(A, \frac{3\kappa
r}4)^c}G_D(x,y) dy\\
&\ge& c_2\, \kappa^\alpha \,r^{-d}\, \frac{\ell((2r)^{-2})}{\ell((
\frac{\kappa r}{2})^{-2})}\int_{D \cap B(A, \frac{3\kappa
r}4)^c}G_D(x,y) dy,
\end{eqnarray*}
for some positive constants $c_1$ and $c_2$. On the other hand,
applying Theorem \ref{ksv-T:Har} we get
\[
\int_{B(A, \frac{3\kappa r}4)} G_D(x,y)   dy\le c_3 \int_{B(A,
\frac{3\kappa r}4)}  G_D(x,A)dy \,\le\,c_4\,r^{d}\,\kappa^d
G_D(x,A),
\]
for some positive constants $c_3$ and $c_4$. Combining these two
estimates we get that
\begin{equation}\label{ksv-efe1}
\int_{D} G_D(x,y)   dy    \,\le\, c_5\,\left(r^{d}\kappa^d+r^{d}
\kappa^{-\alpha}\frac{\ell((\frac{\kappa r}{2})^{-2})}
{\ell((2r)^{-2})}\right)\, G_D(x,A)
\end{equation}
for some constant $c_5>0$.

Let $\Omega=D\setminus \overline{B(A, \frac{\kappa r}2)}$. Note that for any
$z\in B(A, \frac{\kappa r}4)$ and $y\in \Omega$, $ 2^{-1}|y-z|\le|y-A|\le 2|y-z|$.
Thus we get from (\ref{ksv-PK}) and \eqref{ksv-H:1} that for $z\in B(A,\frac{\kappa r}4)$,
\begin{equation}\label{ksv-e:KK1}
c_6^{-1}K_{\Omega}(x, A)  \,\le \,K_{\Omega}(x, z) \,\le\,
c_6K_{\Omega}(x, A)
\end{equation}
for some $c_6>1$. Using the harmonicity of $G_D(\cdot, A)$ in
$D\setminus\{A\}$ with respect to $X$, we can split $G_D(\cdot, A)$
into two parts:
\begin{eqnarray*}
\lefteqn{G_D(x, A)
={\mathbb E}_x \left[G_D(X_{\tau_{\Omega}},A)\right]}\\
&=&{\mathbb E}_x \left[G_D(X_{\tau_{\Omega}},A):\,X_{\tau_{\Omega}} \in B(A,
\frac{\kappa r}4)  \right]\\
& & +
{\mathbb E}_x\left[G_D(X_{\tau_{\Omega}},A):\,X_{\tau_{\Omega}}
\in \{\frac{\kappa r}4\le |y-A|\le \frac{\kappa r}2\}\right]\\
& := &I_1+I_2.
\end{eqnarray*}
Since $G_D(y,A)\le G(y,A)$, by using (\ref{ksv-e:KK1})  and Theorem \ref{ksv-t:Gorigin}, we have
\begin{eqnarray*}
I_1 &\le &
c_6\,K_{\Omega}(x,A) \int_{B(A, \frac{\kappa r}4)}G_D(y,
A)dy \\
& \le & c_7 \,K_{\Omega}(x,A) \int_{B(A, \frac{\kappa r}4)}
\frac{1}{|y-A|^{d-\alpha}\ell(|y-A|^{-2})}\, dy\, ,
\end{eqnarray*}
for some constant $c_7>0$. Since $|y-A|\le 4r \le 4 $, by
\eqref{ksv-el1},
\begin{equation}\label{ksv-efe}
\frac{|y-A|^{\alpha/2}}{\ell(|y-A|^{-2})} \,\le\, c_8 \,
\frac{(4r)^{\alpha/2}}{\ell((4r)^{-2})}
\end{equation}
for some constant $c_8>0$. Thus
\begin{eqnarray*}
I_1 &\le & c_7\, c_8\,K_{\Omega}(x,A) \int_{B(A, \frac{\kappa
r}4)}\frac{1}{|y-A|^{d-\alpha/2}}
\frac{(4r)^{\alpha/2}}{\ell((4r)^{-2})}dy\\
& \le &
c_9\kappa^{\alpha/2}r^{\alpha}\frac1{\ell((4r)^{-2})}K_{\Omega}(x, A)
\end{eqnarray*}
for some constant $c_9>0$. Now   using (\ref{ksv-e:KK1}) again, we  get
\begin{eqnarray*}
  I_1      &\le &
c_{10}\kappa^{\alpha/2-d}r^{\alpha-d}\frac1{\ell((4r)^{-2})}\int_{B(A, \frac{\kappa r}4)} K_{\Omega}(x, z)dz\\
&=&c_{10}\kappa^{\alpha/2-d}r^{\alpha-d}\frac1{\ell((4r)^{-2})}{\mathbb P}_x\left(X_{\tau_{\Omega}}\in B(A, \frac{\kappa r}{4}\right)
\end{eqnarray*}
for some constant $c_{10}>0$. On the other hand, again by Theorem \ref{ksv-t:Gorigin} and \eqref{ksv-efe},
\begin{eqnarray*}
I_2 &=& \int_{\{\frac{\kappa r}4\le |y-A|\le \frac{\kappa r}2\}}
G_{D}(y,A) {\mathbb P}_x(X_{\tau_{\Omega}} \in dy)   \\
&\le &
 c_{11}\int_{\{\frac{\kappa r}4\le |y-A|\le \frac{\kappa r}2\}}
\frac1{|y-A|^{d-\alpha}} \,\frac{1}{\ell(|y-A|^{-2})}
{\mathbb P}_x(X_{\tau_{\Omega}} \in dy)\\
&\le &  c_{12}
\kappa^{\alpha/2-d}\,r^{\alpha-d} \, \frac1{\ell((4r)^{-2})}{\mathbb P}_x
\left(X_{\tau_{\Omega}} \in \{\frac{\kappa r}4\le |y-A|\le
\frac{\kappa r}2\}\right),
\end{eqnarray*}
for some constants $c_{11}>0$ and $c_{12}>0$.

Therefore
$$
G_D(x, A) \,\le\, c_{13}\,
\kappa^{\alpha/2-d}\,r^{\alpha-d}\frac1{\ell((4r)^{-2})}\,
{\mathbb P}_x\left(X_{\tau_{\Omega}} \in B(A, \frac{\kappa r}2)\right)
$$
for some constant $c_{13}>0$. Combining the above with \eqref{ksv-efe1},
we get
\begin{eqnarray*}
\lefteqn{\int_{D} G_D(x,y)   dy}\\
&\le & c_{14}\, r^{\alpha} \,\kappa^{-d-\alpha/2}\,
\frac1{\ell((4r)^{-2})}\left(1+ \frac{\ell((\frac{\kappa
r}{2})^{-2})}{\ell((2r)^{-2})}\right){\mathbb P}_x \left(X_{\tau_{D\setminus
B(A, \frac{\kappa r}2)}} \in B(A, \frac{\kappa r}2)\right),
\end{eqnarray*}
for some constant $c_{14}>0$. It follows immediately that
\begin{eqnarray*}
\lefteqn{\int_{D} G_D(x,y)   dy }\\
& \le & c_{14}\, r^{\alpha} \,
\kappa^{-d-\alpha/2}\,\frac1{\ell((4r)^{-2})}\left(1+
\frac{\ell((\frac{\kappa r}{2})^{-2})} {\ell((2r)^{-2})} \right)
{\mathbb P}_x\left(X_{\tau_{D\setminus B(A, \kappa r)}} \in B(A, \kappa
r)\right).
\end{eqnarray*}
\qed

Combining Lemmas \ref{ksv-l2.1}-\ref{ksv-l2.1_1} and using the translation
invariant property, we have the following

\begin{lemma}\label{ksv-l2.3}
There exists $C_{21}>0$ such that for any open set $D$ with $B(A, \kappa r)\subset D\subset B(Q, r)$
for some $r\in(0, 1)$ and
$\kappa\in (0, 1)$, we have that for every $ x\in D\cap B(Q, \frac{r}2)$,
\begin{eqnarray*}
\lefteqn{{\mathbb P}_x\left(X_{\tau_{D}} \in B(Q, r)^c\right)}\\
& \le &
C_{21}\,\kappa^{-d-\alpha/2 }\, \frac{\ell(r^{-2})}{ \ell((4r)^{-2})}\,
\left(1+\frac{\ell((\frac{\kappa r}{2})^{-2})}{\ell((2r)^{-2})}
\right)  {\mathbb P}_x\left(X_{\tau_{D\setminus B(A, \kappa r) }} \in B(A,
\kappa r) \right).
\end{eqnarray*}
\end{lemma}

Let  $A(x, a,b):=\{ y \in {\mathbb R}^d: a \le |y-x| <b  \}.$

\begin{lemma}\label{ksv-l2.U}
Let $D$ be an open set and $r\in (0,1/2)$. For every $Q \in {\mathbb R}^d$ and
any positive function $u$ vanishing on  $D^c \cap B(Q,
\frac{11}6r)$,
 there is a $\sigma\in
(\frac{10}{6}r, \frac{11}{6}r)$ such that for any $ x \in D
 \cap B(Q, \frac{3}{2}  r)$,
\begin{equation}\label{ksv-e:l2.U}
{\mathbb E}_x\left[u(X_{\tau_{D \cap B(Q, \sigma)}}); X_{\tau_{D \cap B(Q,
\sigma)}} \in B(Q, \sigma)^c\right]
 \le C_{22}\,\frac{r^{\alpha}}{\ell((2r)^{-2})}
\int_{B(Q, \frac{10r}6)^c} J(y-Q)u(y)dy
\end{equation}
for some constant $C_{22}>0$ independent of $Q$ and $u$.
\end{lemma}

\pf Without loss of generality, we may assume that $Q=0$. Note that
by \eqref{ksv-el6}
\begin{eqnarray*}
&&\int^{\frac{11}{6}r}_{\frac{10}{6}r}\int_{A(0, \sigma, 2r)}
\ell((|y|- \sigma)^{-2})^{1/2} (|y|-\sigma)^{-{\alpha}/2}
u(y)\, dy \, d\sigma\\
 &&=\int_{A(0, \frac{10}{6}r , 2r)}
\int^{ |y| \wedge \frac{11}{6}r}_{\frac{10}{6}r}\ell((|y|-
\sigma)^{-2})^{1/2}
(|y|-\sigma)^{-{\alpha}/2}\, d\sigma\, u(y )\,dy \\
&& \le  \int_{A(0, \frac{10}{6}r , 2r)} \left(\int^{ |y| -
\frac{10}{6}r}_{0}\ell(s^{-2})^{1/2}
s^{-{\alpha}/2}ds \right)u(y)dy \\
&&\le c_1 \int_{A(0, \frac{10r}6, 2r)} \ell\left(\left(|y|-  \frac{10r}6
\right)^{-2}\right)^{1/2} \left(|y|- \frac{10r}6\right)^{1-{\alpha}/2} u(y)dy
\end{eqnarray*}
for some positive constant $c_1$.
Using
\eqref{ksv-el2} and \eqref{ksv-el7},
we get that there are constants $c_2>0$ and $c_3>0$ such that
\begin{eqnarray*}
\lefteqn{\int_{A(0, \frac{10r}6, 2r)}\ell\left(\left(|y|-  \frac{10r}6 \right)^{-2}\right)^{1/2}
\left(|y|- \frac{10r}6\right)^{1-{\alpha}/2} u(y)dy }\\
&\le & c_3 \int_{A(0,
\frac{10r}6, 2r)} \ell(|y|^{-2})^{1/2} |y|^{1-{\alpha}/2} u(y)dy\\
&\le & c_3 \frac{r^{1-\alpha/2}}{\ell((2r)^{-2})^{1/2}}  \int_{A(0,
\frac{10r}6, 2r)} \ell(|y|^{-2}) u(y)dy\, .
\end{eqnarray*}
Thus, by taking $c_4>6 c_1 c_3$,
we can conclude that there is a $\sigma\in (\frac{10}{6}r,
\frac{11}{6}r)$ such that
\begin{eqnarray}\label{ksv-e:int}
\lefteqn{\int_{A(0, \sigma, 2r)}\ell((|y|- \sigma)^{-2})^{1/2}\,
(|y|-\sigma)^{-{\alpha}/2}u(y)dy}\nonumber\\
&\le &
c_4\,\frac{r^{-\alpha/2}}{\ell((2r)^{-2})^{1/2}} \int_{A(0,
\frac{10r}6, 2r)} \ell(|y|^{-2}) u(y)dy.
\end{eqnarray}

Let  $x \in D \cap B(0,  \frac{3}{2}  r)$. Note that, since $X$
satisfies the hypothesis ${\bf H}$ in \cite{Sz1}, by Theorem 1 in
\cite{Sz1} we have
\begin{eqnarray*}
&& {\mathbb E}_x\left[u(X_{\tau_{D \cap B(0, \sigma)}}); X_{\tau_{D \cap B(0,
\sigma)}} \in
B(0, \sigma)^c \right]\\
 &&= {\mathbb E}_x\left[u(X_{\tau_{D \cap B(0, \sigma)}}); X_{\tau_{D \cap
B(0, \sigma)}} \in B(0, \sigma)^c, \, \tau_{D \cap
 B(0, \sigma)} =\tau_{B(0, \sigma)}  \right]\\
&&= {\mathbb E}_x\left[u(X_{\tau_{ B(0, \sigma)}}); X_{ \tau_{B(0, \sigma)}}
\in
B(0, \sigma)^c, \, \tau_{D \cap B(0, \sigma)} =
 \tau_{B(0, \sigma)}  \right]\\
&&\le {\mathbb E}_x\left[u(X_{\tau_{ B(0, \sigma)}}); X_{\tau_{B(0, \sigma)}}
\in B(0, \sigma)^c  \right] \,=\,\int_{B(0, \sigma)^c}K_{B(0,
\sigma)}(x, y)u(y)dy.
\end{eqnarray*}
Since  $ \sigma  <2r < 1$, from
\eqref{ksv-P1} in Proposition \ref{ksv-p:Poisson1}, Proposition \ref{ksv-p:Poisson3} we have
\begin{eqnarray*}
\lefteqn{ {\mathbb E}_x\left[u(X_{\tau_{D \cap B(0, \sigma)}}); X_{\tau_{D \cap B(0,
\sigma)}} \in B(0, \sigma)^c \right]
 \,\le\,
\int_{  B(0, \sigma)^c   }
 K_{B(0, \sigma)}(x, y)u(y)dy}\\
&\le &\, c_5 \int_{A(0, \sigma, 2r)}
\frac{\sigma^{\alpha/2-d}}{\left(\ell(\sigma^{-2})\right)^{1/2}}
\frac{(\ell((|y|-\sigma)^{-2}))^{1/2}} {( |y|-\sigma)^{\alpha/2}}
u(y)dy\\
&& + c_5 \int_{B(0, 2r)^c}j(|y|-\sigma)
\frac{\sigma^{\alpha/2}}
 {(\ell(\sigma^{-2}))^{1/2}}\frac{(\sigma-|x|)^{\alpha/2}}
{(\ell((\sigma-|x|)^{-2}))^{1/2}}
 u(y)dy
\end{eqnarray*}
for some constant $c_5>0$.

When $y \in A(0,  2r , 4)$ we have $\frac1{12}|y|\le |y|-\sigma $, while when $|y|\ge 4$ we have
$|y|-\sigma\ge |y|-1$. Since $ \sigma-|x|\le\sigma \le {2r}$, we
have by \eqref{ksv-el1} and the monotonicity of $j$,
$$
j(|y|-\sigma) \frac{\sigma^{\alpha/2}}{(\ell(\sigma^{-2}))^{1/2}}
\frac{(\sigma-|x|)^{\alpha/2}} {(\ell((\sigma-|x|)^{-2}))^{1/2}}
 \,\le\, c_6 j\left(\frac{|y|}{12}\right)
\frac{r^{\alpha}}{\ell((2r)^{-2})} , \quad y \in A(0,  2r , 4)
$$
and
$$
j(|y|-\sigma) \frac{\sigma^{\alpha/2}}{(\ell(\sigma^{-2}))^{1/2}}
\frac{(\sigma-|x|)^{\alpha/2}} {(\ell((\sigma-|x|)^{-2}))^{1/2}}
 \,\le\, c_6 j(|y|-1)
\frac{r^{\alpha}}{\ell((2r)^{-2})} , \quad |y|\ge 4
$$
for some constant $c_6>0$. Thus by applying \eqref{ksv-H:1} and
\eqref{ksv-H:2}, we get
$$
j(|y|-\sigma) \frac{\sigma^{\alpha/2}}{(\ell(\sigma^{-2}))^{1/2}}
\frac{(\sigma-|x|)^{\alpha/2}} {(\ell((\sigma-|x|)^{-2}))^{1/2}}
 \,\le\, c_7    j(|y|)
\frac{r^{\alpha}}{\ell((2r)^{-2})}
$$
for some constant $c_7>0$.
Therefore,
\begin{eqnarray*}
\lefteqn{\int_{B(0, 2r)^c}j(|y|-\sigma)\frac{\sigma^{\alpha/2}} {(\ell(\sigma^{-2}))^{1/2}}\frac{(\sigma-|x|)^{\alpha/2}}
{(\ell((\sigma-|x|)^{-2}))^{1/2}}
 u(y)dy}\\
 & \le &  c_5 c_7 \frac{r^{\alpha}}{\ell((2r)^{-2})} \int_{B(0, 2r)^c} J(y) u(y)\, dy\, .
\end{eqnarray*}

On the other hand, by \eqref{ksv-el1}, \eqref{ksv-e:int}  and Theorem \ref{ksv-t:Jorigin}, we have that
\begin{eqnarray*}
\lefteqn{\int_{A(0, \sigma, 2r)}
\frac{\sigma^{\alpha/2-d}}{\left(\ell(\sigma^{-2})\right)^{1/2}}
\frac{(\ell((|y|-\sigma)^{-2}))^{1/2}}
{( |y|-\sigma)^{\alpha/2}} u(y)dy}\\
&\le& \left(\frac{10r}{6}\right)^{-d}\frac{\sigma^{\alpha/2}}
{\left(\ell(\sigma^{-2})\right)^{1/2}} \int_{A(0, \sigma, 2r)}
   \frac{(\ell((|y|-\sigma)^{-2}))^{1/2}}
{( |y|-\sigma)^{\alpha/2}} u(y)dy\\
&\le & c_8
r^{-d}\frac{(2r)^{\alpha/2}}{\left(\ell((2r)^{-2})\right)^{1/2}}
\,\frac{r^{-\alpha/2}}{\left(\ell((2r)^{-2})\right)^{1/2}}
\int_{A(0, \frac{10r}6, 2r)} \ell(|y|^{-2}) u(y)dy \\
&\le & c_{9} \frac{r^{\alpha}}{\ell((2r)^{-2})} \int_{A(0,
\frac{10r}6, 2r)}\ell(|y|^{-2}) |y|^{-d-\alpha} u(y) dy\\
&\le &  c_{10} \frac{r^{\alpha}}{\ell((2r)^{-2})}
\int_{A(0, \frac{10r}6, 2r)}J(y) u(y) dy
\end{eqnarray*}
for some positive constants $c_8$, $c_9$ and $c_{10}$.
Hence, by combining the last two displays we arrive at
$$
{\mathbb E}_x\left[u(X_{\tau_{D \cap B(0, \sigma)}}); X_{\tau_{D \cap B(0,
\sigma)}} \in
B(0, \sigma)^c \right]\\
 \,\le\,
c_{11}\,\frac{r^{\alpha}}{\ell((2r)^{-2})}
\int_{B(0, \frac{10r}6)^c}J(y)u(y)dy
$$
for some constant $c_{11}>0$.
 \qed

\begin{lemma}\label{ksv-l2.2}
Let $D$ be an open set and $r\in (0,1/2)$. Assume that $B(A, \kappa r)\subset D\cap
B(Q, r)$ for $\kappa\in (0, 1/2 ]$.
Suppose that $u\ge0$ is regular harmonic in $D\cap B(Q, 2r)$ with respect to
$X$ and $u=0$ in $D^c\cap B(Q, 2r)$. If $w$ is a regular harmonic
function with respect to $X$ in $D\cap B(Q, r)$ such that
$$
w(x)=\left\{
\begin{array}{ll}
u(x), & x\in B(Q, \frac{3r}2)^c\cup (D^c\cap B(Q, r)),\\
0, & x \in A(Q, r, \frac{3r}2),
\end{array}\right.
$$
then
$$
u(A) \ge w(A) \ge C_{23}\,\kappa^{\alpha} \frac{\ell((2r)^{-2})}
{\ell((\kappa r)^{-2})} \,u(x), \quad x \in D \cap B(Q,\frac32 r)
$$
for some constant $C_{23}>0$.
\end{lemma}

\pf
Without loss of generality we may assume $Q=0$. Let  $x \in D \cap B(0,\frac32 r)$.
The left hand side inequality in the conclusion of the lemma is clear from the fact that $u$ dominates $w$ on $(D\cap B(0,r))^c$ and both functions are regular harmonic in $D\cap B(0,r)$.
Thus we only need to prove the right hand side inequality.
By Lemma \ref{ksv-l2.U} there exists $\sigma\in (\frac{10r}6, \frac{11r}6)$ such that \eqref{ksv-e:l2.U} holds.
Since $u$ is regular harmonic in $D\cap B(0, 2r)$ with respect to $X$ and equal to zero on $D^c\cap B(0,2r)$,
it follows that
\begin{equation}\label{ksv-e:l2.2}
u(x)= {\mathbb E}_x\left[u(X_{\tau_{D \cap B(0, \sigma)}}); \,X_{\tau_{D \cap
B(0, \sigma)}} \in B(0, \sigma)^c \right] \le
c_1\frac{r^{\alpha}}{\ell((2r)^{-2})} \int_{B(0, \frac{10r}6)^c} J(y)u(y)dy
\end{equation}
for some constant $c_1>0$.
On the other hand, by \eqref{ksv-P2} in Proposition \ref{ksv-p:Poisson1},
we have  that
\begin{eqnarray*}
w(A)&=& \int_{B(0, \frac{3r}2)^c}
 K_{D\cap B(0, r)}(A, y)u(y)dy
\ge  \int_{B(0, \frac{3r}2)^c}
 K_{B(A, \kappa r)}(A, y)u(y)dy\\
&\ge &  c_2 \int_{B(0, \frac{3r}2)^c} J(A-y) \frac{(\kappa
r)^\alpha}{\ell((\kappa r)^{-2})}  u(y)dy
\end{eqnarray*}
for some constant $c_2>0$. Note that $|y-A|\le  2|y|$ in $
A(0,\frac{3r}2, 4)$ and that $|y-A|\le |y|+1$ for $|y|\ge 4$. Hence
by the monotonicity of $j$, \eqref{ksv-H:1} and \eqref{ksv-H:2},
$$
w(A)\,\ge\, c_3\,\frac{(\kappa r)^\alpha}{\ell((\kappa r)^{-2})}
\int_{B(0, \frac{3r}2)^c} J( y)
 u(y) dy
$$
for some constant $c_3>0$. Therefore, by \eqref{ksv-e:l2.2}
$$
w(A)\ge c_4 c_1^{-1} \,\kappa^{\alpha}  \frac{\ell((2r)^{-2})}
{\ell((\kappa r)^{-2})}   \,u(x)\, .
$$
\qed

\begin{defn}\label{ksv-fat}
Let $\kappa \in (0,1/2]$. We say that an open set $D$ in ${\mathbb R}^d$ is
$\kappa$-fat if there exists $R>0$ such that for each $Q \in
\partial D$ and $r \in (0, R)$, $D \cap B(Q,r)$ contains a ball
$B(A_r(Q),\kappa r)$. The pair $(R, \kappa)$ is called the
characteristics of the $\kappa$-fat open set $D$.
\end{defn}

Note that all Lipschitz domain and all non-tangentially accessible
domain (see \cite{JK} for the definition) are $\kappa$-fat. The
boundary of a $\kappa$-fat open set can be highly nonrectifiable
and, in general, no regularity of its boundary can be inferred.
Bounded $\kappa$-fat open set may be disconnected.

Since $\ell$ is slowly varying at $\infty$, we get the Carleson's
estimate from Lemma \ref{ksv-l2.2}.

\begin{corollary}\label{ksv-c:Carl}
Suppose that $D$ is  a $\kappa$-fat open set with the
characteristics $(R, \kappa)$. There exists a constant $C_{24}$
depending on the characteristics $(R,\kappa)$ such that if $r \le
R\wedge \frac12$, $Q\in \partial D$, $u\ge0$ is regular harmonic in
$D\cap B(Q, 2r)$ with respect to $X$ and $u=0$ in $D^c\cap B(Q,
2r)$, then
$$
u\left(A_r(Q)\right) \,\ge C_{24}\, u(x)\, , \quad  \forall
x \in D \cap B(Q,\frac32 r)\, .
$$
\end{corollary}

The next theorem is a boundary Harnack principle for (possibly
unbounded) $\kappa$-fat open set and it is the main result of this
subsection.

\begin{thm}\label{ksv-BHP}
Suppose that $D$ is  a $\kappa$-fat open set with the
characteristics $(R, \kappa)$. There exists a  constant $C_{25}>1$
depending on the characteristics $(R,\kappa)$ such that if $r \le
R\wedge \frac14$ and $Q\in\partial D$, then for any nonnegative
functions $u, v$ in ${\mathbb R}^d$ which are regular harmonic in $D\cap B(Q,
2r)$ with respect to $X$ and vanish in $D^c \cap B(Q, 2r)$, we have
$$
C_{25}^{-1}\,\frac{u(A_r(Q))}{v(A_r(Q))}\,\le\, \frac{u(x)}{v(x)}\,\le
C_{25}\,\frac{u(A_r(Q))}{v(A_r(Q))}, \qquad x\in D\cap B(Q, \frac{r}2)\, .
$$
\end{thm}

\pf
Since $\ell$ is slowly varying at $\infty$ and
locally bounded above and below by positive constants,
there exists a
constant $c>0$ such that for every $r\in (0,1/4)$,
\begin{equation}\label{ksv-lll}
\max \left(\frac{\ell(r^{-2})}{ \ell((\kappa r)^{-2})    },\,
\frac{\ell((2r)^{-2})}{ \ell((4r)^{-2})},\, \frac{\ell((\frac{\kappa
r}{2})^{-2})}{\ell((4r)^{-2})},\, \frac{\ell((\kappa
r)^{-2})}{\ell((2r)^{-2})} \right) \,\le\, c\, .
\end{equation}

Fix $r\in (0, R\wedge \frac14)$ throughout this proof. Without loss of generality
we may assume that $Q=0$ and $u(A_r(0))=v(A_r(0))$. For simplicity,
we will write $A_r(0)$ as $A$ in the remainder of this proof. Define
$u_1$ and $u_2$ to be regular harmonic functions in $D\cap B(0, r)$
with respect to $X$ such that
$$
u_1(x)=\left\{
\begin{array}{ll}
u(x), & x\in A(0,r, \frac{3r}{2}),\\
0, & x\in B(0, \frac{3r}2)^c\cup(D^c\cap B(0, r))
\end{array}
\right.
$$
and
$$
u_2(x)=\left\{
\begin{array}{ll}
0, &  x\in A(0,r, \frac{3r}{2}), \\
u(x), & x\in B(0, \frac{3r}2)^c\cup(D^c\cap B(0, r)).
\end{array}
\right.
$$
and note that $u=u_1+u_2$. If $D\cap A(0,r,\frac{3r}{2})=\emptyset$,
then $u_1=0$ and the inequality (\ref{ksv-e2.6}) below holds
trivially. So we assume that $D\cap A(0,r,\frac{3r}{2})$ is not empty.
Then by Lemma \ref{ksv-l2.2},
$$
u(y)\le c_1 \kappa^{-\alpha} \frac{\ell((\kappa
r)^{-2})}{\ell((2r)^{-2})}\, u(A), \qquad y\in D\cap B(0, \frac{3r}2),
$$
for some constant $c_1>0$.
For $x\in D\cap B(0, \frac{r}2)$, we have
\begin{eqnarray*}
u_1(x)&=& {\mathbb E}_x\left[u(X_{\tau_{D\cap B(0, r)}}): X_{\tau_{D\cap B(0,
r)}}\in
D\cap A(0,r,\frac{3r}{2})\right]\\
&\le&\left(\sup_{D\cap A(0,r,\frac{3r}{2})}u(y)\right)
   {\mathbb P}_x\left( X_{\tau_{D\cap B(0, r)}}\in
D\cap A(0,r,\frac{3r}{2})\right)              \\
&\le&\left(\sup_{D\cap A(0,r,\frac{3r}{2})}u(y)\right)
   {\mathbb P}_x\left( X_{\tau_{D\cap B(0, r)}}\in
B(0,r)^c        \right)     \\
&\le&c_1\,\kappa^{-\alpha} \frac{\ell((\kappa
r)^{-2})}{\ell((2r)^{-2})} \,u(A) \,{\mathbb P}_x\left( X_{\tau_{D\cap B(0,
r)}}\in B(0,r)^c \right).
\end{eqnarray*}
Now using Lemma \ref{ksv-l2.3} (with $D$ replaced by $D\cap B(0,r)$) and \eqref{ksv-lll}, we have that for $ x\in D\cap B(0, \frac{r}2)$,
\begin{eqnarray}
&&u_1(x)\\
&&\le\, c_2\,\kappa^{-d-\frac32\alpha }\, \frac {\ell((\kappa
r)^{-2})}{\ell((2r)^{-2})}\frac{\ell(r^{-2})}{\ell((4r)^{-2})}\,
\left(1+\frac{\ell((\frac{\kappa r}{2})^{-2})}{\ell((4r)^{-2})}
\right)\,u(A)\times \nonumber\\
&&\quad \times\ {\mathbb P}_x\left( X_{\tau_{(D\cap B(0,r))\setminus
B(A, \frac{\kappa r}2)}} \in B(A, \frac{\kappa r}2)\right) \nonumber\\
&&\le\,c_3 \,u(A)\,{\mathbb P}_x\left( X_{\tau_{(D\cap B(0,r))\setminus B(A,
\frac{\kappa r}2)}} \in B(A, \frac{\kappa r}2)\right) \label{ksv-e2.3}
\end{eqnarray}
for some positive constants $c_2$ and $c_3=c_3(\kappa)$.
Since $r <1/4$, Theorem \ref{ksv-T:Har} implies that
$$
u(y)\,\ge\, c_4\,u(A), \qquad y\in B(A, \frac{\kappa r}2)
$$
for some constant $c_4>0$. Therefore for $x\in D\cap B(0,
\frac{r}2)$
\begin{equation}\label{ksv-e2.4}
u(x) \,=\,  {\mathbb E}_x\left[u(X_{\tau_{(D\cap B(0, r))\setminus B(A,
\frac{\kappa r}2)}}) \right]     \,\ge\, c_4\,u(A)\,
{\mathbb P}_x\left(X_{\tau_{(D\cap B(0,r))\setminus B(A, \frac{\kappa r}2)}}
\in B(A, \frac{\kappa r}2)\right).
\end{equation}
Using (\ref{ksv-e2.3}), the analogue of (\ref{ksv-e2.4}) for $v$, and the
assumption that $u(A)=v(A)$, we get that for $x\in D\cap B(0,
\frac{r}2)$,
\begin{equation}\label{ksv-e2.6}
u_1(x)\,\le \,c_3\,v(A)\, {\mathbb P}_x\left(X_{\tau_{(D\cap B(0, r))
\setminus B(A, \frac{\kappa r}2)}} \in B(A, \frac{\kappa
r}2)\right)\,\le \,c_5\,v(x)
\end{equation}
for some constant $c_5=c_5(\kappa)>0.$ For $x\in D\cap B(0, r)$, we
have
\begin{eqnarray*}
u_2(x)&=& \int_{B(0, \frac{3r}2)^c}K_{D\cap B(0, r)}
 (x, z)u(z)dz\\
&=& \int_{B(0, \frac{3r}2)^c}
 \int_{D\cap B(0, r)} G_{D\cap B(0, r)}(x, y)  J(y-z)dy\, u(z)\, dz.
\end{eqnarray*}
Let
$$
s(x):=\int_{D\cap B(0, r)}G_{D\cap B(0, r)}(x, y)dy.
$$
Note that for every $y \in B(0,r)$ and $z \in B(0, \frac{3r}2)^c$,
$$
\frac13|z| \,\le\, |z|-r \,\le\, |z|-|y| \, \le\, |y-z|\, \le \,
|y|+|z|\, \le \, r+|z| \le 2 |z|\, ,
$$
and that for every $y \in B(0,r)$ and $z \in B(0, 12)^c$,
$$
|z|-1\, \le \, |y-z|\, \le \,|z|+1.
$$
So by the monotonicity of $j$, for every $y \in B(0,r)$ and $z \in
A(0, \frac{3r}2, 12)$,
$$
j(12|z|) \, \le \,j(2|z|) \, \le \,J(y-z) \, \le \, j\left(\frac{|z|}{3}\right) \,
\le \,j\left(\frac{|z|}{12}\right)\, ,
$$
and for every $y \in B(0,r)$ and every $z \in B(0,12)^c$,
$$
j(|z|-1)\, \le \,J(y-z) \, \le \, j(|z|+1).
$$
Using \eqref{ksv-H:1} and \eqref{ksv-H:2}, we have that, for every $y \in
B(0,r)$ and $z \in B(0, \frac{3r}2)^c$,
$$  c_6^{-1} j(|z|) \, \le  \,J(y-z) \, \le \,c_6\,j(|z|)
$$
for some constant $c_6>0$. Thus we have
\begin{equation}\label{ksv-e2.5}
c_7^{-1}\le \left(\frac{u_2(x)}{u_2(A)}\right)\left(\frac{s(x)}{s(A)}\right)^{-1}\le c_7,
\end{equation}
for some constant $c_7>1$.
 Applying (\ref{ksv-e2.5}) to $u$, and $v$ and Lemma \ref{ksv-l2.2} to $v$ and
$v_2$, we obtain for $x\in D\cap B(0, \frac{r}2)$,
\begin{eqnarray}\label{ksv-e2.7}
u_2(x) &\le & c_7\,u_2(A)\,\frac{s(x)}{s(A)}\,\le\, c_{7}^2\,
\frac{u_2(A)}{v_2(A)}\,v_2(x)\, \le\, c_{8}\, \kappa^{-\alpha}
\frac{\ell((\kappa r)^{-2})}
{\ell((2r)^{-2})}\frac{u(A)}{v(A)}\,v_2(x)\nonumber \\
& = &
c_{8}\,\kappa^{-\alpha} \frac{\ell((\kappa
r)^{-2})}{\ell((2r)^{-2})}\,v_2(x),
\end{eqnarray}
for some constant $c_8>0.$
 Combining (\ref{ksv-e2.6}) and (\ref{ksv-e2.7}) and applying \eqref{ksv-lll}, we
have
$$
u(x)\,\le\, c_{9} \,v(x), \qquad x\in D\cap B(0, \frac{r}2),
$$
for some constant $c_{9}=c_{9}(\kappa)>0.$ \qed

\medskip
\noindent {\bf Acknowledgment:} We thank Qiang Zeng for his comments
on the first version of this paper.
We also thank the referee for helpful comments.

\vspace{.5in}
\begin{singlespace}

\end{singlespace}
\end{doublespace}


\small
\begin{thebibliography}{99}








\bibitem{BL02a} Bass, R.F., Levin, D.A.:
Harnack inequalities for jump processes. {\em Potential Anal.}, {\bf
17} (2002), 375--388.



\bibitem{Ber} Bertoin, J.: {\em L\'evy Processes}. Cambridge University
Press, Cambridge, 1996.



\bibitem{BGT} Bingham, N.H., Goldie, C.M., Teugels, J.L.:
{\em Regular Variation}. Cambridge University Press, Cambridge,
1987.



\bibitem{Bog97} Bogdan, K.:
The boundary Harnack principle for the fractional Laplacian. {\em
Studia Math.}, {\bf 123} (1997), 43--80.

\bibitem{B3} Bogdan, K.: Sharp estimates for the Green function
in Lipschitz domains. {\em J. Math. Anal. Appl. \bf 243} (2000),
326-337.

\bibitem{BBKRSV} Bogdan, K.,   Byczkowski, T.,  Kulczycki, T.,
Ryznar, M.,  Song, R. and Vondra\v{c}ek, Z.: {\em Potential analysis
of stable processes and its extensions}.  Lecture Notes in
Mathematics, {\bf 1980}. Springer-Verlag, Berlin, 2009.



\bibitem{C} Chen, Z.-Q.: On notions of harmonicity.  {\em Proc. Amer. Math. Soc.}
{\bf  137} (2009),  3497--3510.


\bibitem{CKS1} Chen, Z.-Q., Kim, P.  and Song, R.:
Heat kernel estimates for Dirichlet fractional Laplacian. {\em J.
European Math. Soc.}, {\bf 12} (2010), 1307--1329.

\bibitem{CKS2} Chen, Z.-Q., Kim, P.  and Song, R.:
Sharp heat kernel estimates for relativistic stable processes in
open sets. {\em Ann. Probab.}, to appear, 2011.

\bibitem{CKS3} Chen, Z.-Q., Kim, P.  and Song, R.:
Dirichlet heat kernel estimates for $\Delta^{\alpha/2}
+\Delta^{\beta/2}$. {\em Ill J. Math.}, to appear, 2011.

\bibitem{CKS4} Chen, Z.-Q., Kim, P.  and Song, R.:
Heat kernel estimates for $\Delta + \Delta^{\alpha/2}$ in $C^{1, 1}$ open sets.
{\em J.
London Math. Soc.}, to appear, 2011.


\bibitem{CKS5} Chen, Z.-Q., Kim, P.  and Song, R.:
Global Heat Kernel Estimates for $\Delta+\Delta^{\alpha/2}$  in Half-space-like domains.
Preprint, 2011.



\bibitem{CKSV} Chen, Z.-Q., Kim, P., Song, R. and Vondra{\v{c}}ek,
Z.: Boundary Harnack principle for $\Delta + \Delta^{\alpha/2}$.
{\em Trans.~Amer.~Math.~Soc.}, to appear 2011.

\bibitem{CKSV2} Chen, Z.-Q., Kim, P., Song, R. and Vondra{\v{c}}ek,
Z.: Sharp Green function estimates for $\Delta+ \Delta^{\alpha/2}$
in $C^{1,1}$ open sets and their applications. {\em Ill. J. Math.},
to appear, 2011.

\bibitem{CK1} Chen, Z.-Q. and  Kumagai, T.:
Heat kernel estimates for stable-like
processes on $d$-sets.  {\it Stoch. Proc. Appl. \bf 108} (2003), 27--62.



\bibitem{CK2}  Chen, Z.-Q. and  Kumagai, T.:
Heat kernel estimates for jump processes of mixed types on metric
measure spaces. {\it Probab. Theory Relat. Fields \bf 140} (2008),
277--317.



\bibitem{CS05} Chen, Z.-Q., Song, R.: Two-sided eigenvalue estimates
for subordinate processes in domains. {\em J. Funct. Anal.}, {\bf
226} (2005), 90--113.

\bibitem{CS06a} Chen, Z.-Q., Song, R.: Continuity of eigenvalues
of subordinate processes in domains. {\em Math. Z.},  {\bf 252}
(2006), 71--89.

\bibitem{CS06c} Chen, Z.-Q., Song, R.: Spectral properties of
subordinate processes in domains. In Stochastic Analysis and Partial
Differential Equations, AMS, Providence, 2007.






\bibitem{Fris} Fristedt, B. E.:
Sample functions of stochastic processes with stationary,
independent increments. {\em  Advances in probability and related
topics}, Vol. 3, pp. 241--396, Dekker, New York, 1974.



\bibitem{GPRSSV} Glover, J., Pop-Stojanovic, Z., Rao, M., \v{S}iki\'{c}, H.,
Song, R., Vondra\v{c}ek, Z.: Harmonic functions of subordinate
killed Brownian motions. {\em J. Funct. Anal.},  {\bf 215} (2004),
399--426.

\bibitem {GRSS} Glover, J., Rao, M., \v{S}iki\'{c}, H., Song, R.:
$\Gamma$-potentials. In Classical and modern potential theory and
applications (Chateau de Bonas, 1993), 217--232, Kluwer Acad. Publ.,
Dordrecht, 1994.



\bibitem{H}
Hansen, W.:  Uniform boundary Harnack principle and generalized
triangle property. {\em J. Funct. Anal.} {\bf 226}, (2005),
452--484.



\bibitem{JK} Jerison, D. S. and Kenig, C. E.: Boundary behavior
of harmonic functions in non-tangentially accessible domains. {\it
Adv. Math.}, {\bf 46} (1982), 80--147.

\bibitem{KSV1} Kim, P. Song, R. and Vondra\v{c}ek, Z.:
Boundary Harnack principle for subordinate Brownian motion, {\em
Stoch. Proc. Appl.} {\bf 119} (2009), 1601--1631.

\bibitem{KSV2} Kim, P. Song, R. and Vondra\v{c}ek, Z.: On the
potential theory of one-dimensional subordinate Brownian motions
with continuous components. {\it Potential Anal.}, {\bf 33} (2010),
153--173.

\bibitem{KSV3} Kim, P. Song, R. and Vondra\v{c}ek, Z.: Two-sided
Green function estimates for the killed subordinate Brownian
motions. Preprint, 2011.




\bibitem{Kw}
 Kwa\'snicki, M.:
 Spectral analysis of subordinate Brownian motions in half-line.
Preprint, 2010.

\bibitem{Ky} Kyprianou, A. E.:
{\em Introductory lectures on fluctuations of L\'evy processes with
applications}, Springer, Berlin, 2006.



\bibitem{PS71} Port, S.C., Stone, J.C.: Infinitely divisible process
and their potential theory I, II. {\em Ann. Inst. Fourier}, {\bf 21}
(1971), 157--275 and  179--265.




\bibitem{RSV} Rao, M., Song, R., Vondra\v{c}ek, Z.:
Green function estimates and Harnack inequality for subordinate
Brownian motions. {\em Potential Anal.}, {\bf 25} (2006), 1--27.


\bibitem{Ryz02} Ryznar, M.:
Estimates of Green functions for relativistic $\alpha$-stable
process. {\em Potential Anal.},  {\bf 17} (2002), 1-23.

\bibitem{Sat} Sato, K.-I.: {\em L\'evy Processes and Infinitely Divisible
Distributions}. Cambridge University Press, 1999.

\bibitem{SSV} Schilling, R. L.  Song, R. and Vondra{\v{c}}ek, Z.:
{\em Bernstein Functions: Theory and Applications}. de Gruyter
Studies in Mathematics 37. Berlin: Walter de Gruyter, 2010.


\bibitem{SiSV} \v{S}iki\'{c}, H., Song, R., Vondra\v{c}ek, Z.:
Potential theory of geometric stable processes. {\em Probab. Theory
Related Fields},  {\bf 135} (2006), 547--575.

\bibitem{Sil} Silverstein, M. L.: Classification of coharmonic
and coinvariant functions for a L\'evy process. {\em Ann.~Probab.}
{\bf 8} (1980), 539--575,


\bibitem{Son04} Song, R.: Sharp bounds on the density, Green function and
jumping function of subordinate killed BM. {\em Probab. Theory
Related Fields},  {\bf 128} (2004), 606-628.

\bibitem{SV03} Song, R., Vondra\v{c}ek, Z.:
Potential theory of subordinate killed Brownian motion in  a domain.
{\em Probab. Theory Related Fields},  {\bf 125} (2003), 578--592.

\bibitem{SV04} Song, R., Vondra\v{c}ek, Z.:
Harnack inequalities for some classes of Markov processes. {\em
Math. Z.}, {\bf 246} (2004), 177--202.

\bibitem{SV04b} Song, R., Vondra\v{c}ek, Z.: Sharp bounds for Green functions and
jumping functions of subordinate killed Brownian motions in bounded
$C^{1,1}$ domains. {\em Elect. Comm. Probab.},  {\bf 9} (2004),
96--105.



\bibitem{SV06} Song, R., Vondra\v{c}ek, Z.:
Potential theory of special subordinators and subordinate killed
stable processes. {\em J. Theoret. Probab.},  {\bf 19} (2006),
817--847.


\bibitem{SW99} Song, R., Wu, J.-M.:
Boundary Harnack principle for symmetric stable processes. {\em J.
Funct. Anal.},  {\bf 168} (1999), 403--427.

\bibitem{Sz1} Sztonyk, P.:
On harmonic measure for L\'evy processes, {\em Probab. Math.
Statist.}, {\bf 20} (2000), 383--390.

\bibitem{Sz2}
Sztonyk, P.: Boundary potential theory for stable L\'evy processes,
{\em Colloq. Math.}, {\bf 95(2)}  (2003), 191--206.

\bibitem{WYY} Watanabe, S., Yano, K., Yano, Y.: A density formula for
the law of time spent on the positive side of a one-dimensional diffusion process,
\emph{J. Math. Kyoto Univ.} {\bf 45} (2005), 781--806.

\bibitem{Z} Z\"ahle, M.: Potential spaces and traces of L\'evy processes
on $h$-sets, {\em J. Contemp. Math. Anal.} {\bf 44} (2009),
117--145.

\end{thebibliography}
\end{document}